\documentclass[11pt,letterpaper]{amsart}

\usepackage{graphicx}
\usepackage{amssymb}
\usepackage{array}
\usepackage{amsthm}
\usepackage{mathtools}
\usepackage[hidelinks]{hyperref}
\usepackage{tgtermes}
\usepackage{mathrsfs}
\usepackage{stmaryrd}
\usepackage[all]{xy}
\usepackage{verbatim} 
\usepackage{fullpage}  
\usepackage{hyperref}
\usepackage{setspace}
\usepackage{tikz-cd}
\usepackage{cleveref}
\usepackage{amsmath}
\usetikzlibrary{positioning}
\newcommand{\RR}{\mathbb{R}}
\newcommand{\ZZ}{\mathbb{Z}}
\newcommand{\PP}{\mathbb{P}}

\newcommand{\CC}{\mathbb{C}}
\newcommand{\HH}{\mathbb{H}}
\newcommand{\QQ}{\mathbb{Q}}

\newcommand{\Addresses}{{
  \bigskip
  \footnotesize

  \textsc{Department of Mathematics, University of California, Irvine}\par\nopagebreak
  \textit{E-mail address}: \texttt{svraman@uci.edu}
}}

\DeclareMathOperator{\Aut}{Aut}

\DeclareMathOperator{\SO}{SO}
\DeclareMathOperator{\OO}{O}

\DeclareMathOperator{\Hom}{Hom}

\DeclareMathOperator{\Homeo}{Homeo}

\DeclareMathOperator{\rank}{rank}
\DeclareMathOperator{\Gr}{Gr}

\DeclareMathOperator{\END}{End}
\DeclareMathOperator{\MOD}{Mod}

\DeclareMathOperator{\Isom}{Isom}
\DeclareMathOperator{\Vol}{Vol}

\DeclareMathOperator{\inv}{^{-1}}
\DeclareMathOperator{\id}{id}

\DeclareMathOperator{\Teich}{Teich}
\DeclareMathOperator{\Diff}{Diff}
\DeclareMathOperator{\SHomeo}{SHomeo}

\DeclareMathOperator{\SDiff}{SDiff}
\DeclareMathOperator{\LDiff}{LDiff}
\DeclareMathOperator{\SMod}{SMod}
\DeclareMathOperator{\LMod}{LMod}

\DeclareMathOperator{\HK}{HK}
\DeclareMathOperator{\RF}{RF}

\newtheorem{theorem}{Theorem}[section]
\newtheorem{mainthm}{Theorem}

\newtheorem{proposition}[theorem]{Proposition}
\newtheorem{lemma}[theorem]{Lemma}
\newtheorem{corollary}[theorem]{Corollary}
\newtheorem{conjecture}[theorem]{Conjecture}
\theoremstyle{definition}
\newtheorem{definition}[theorem]{Definition}
\newtheorem{question}[theorem]{Question}

\newtheorem{example}[theorem]{Example}
\newtheorem{remark}[theorem]{Remark}

\numberwithin{equation}{subsection}
\title{A smooth Birman--Hilden theory for hyperk\"{a}hler manifolds}
\author{Sidhanth Raman}
\begin{document}
\maketitle
\begin{abstract}
This article grew out of an effort to understand the smooth mapping class groups of certain 4-manifolds in a geometric manner. We prove a smooth analog of the Birman--Hilden theorem for manifolds that admit a hyperk\"{a}hler structure. This allows us to probe the smooth mapping class groups associated to certain manifolds with nontrivial fundamental groups. Along the way, and of independent interest, we prove a global Torelli theorem for generalized Enriques manifolds. The techniques used are analogous to Teichm\"{u}ller-theoretic methods in the classical theory of mapping class groups. We apply this hyperk\"{a}hler Birman--Hilden theorem to obtain results regarding smooth, metric, and complex Nielsen realization on Enriques surfaces.
\end{abstract}

\section{Introduction}

Mapping class groups play a fundamental role in our understanding of algebraic varieties and their moduli. Monodromy in the universal family over the moduli space $\mathscr{M}_X$ gives rise to a natural representation
$$\rho : \pi_1(\mathscr{M}_X)\longrightarrow \MOD(X) := \pi_0(\Diff^+(X)).$$
When $X$ is an algebraic curve, the map $\rho$ is an isomorphism. For higher dimensional varieties, the connection between mapping class groups and moduli spaces is not well understood yet, which can be summarized in the following question:

\begin{question}\label{pi1Mod}
    Let $X$ be the smooth manifold underlying a complex algebraic variety, and $\mathscr{M}_X$ its associated moduli space of complex structures. What is the relationship between $\pi_1(\mathscr{M}_X)$ and $\MOD(X)$? Is $\rho$ injective or surjective?
\end{question}

Unfortunately, our understanding of the global structure of mapping class groups of complex algebraic varieties is still in its infancy; see \cite{hain2023mapping} for related work and open problems on the mapping class groups of simply-connected K\"{a}hler manifolds. The situation for spaces with nontrivial fundamental groups is even more opaque. This article is an attempt to make progress in this direction through the language of covering spaces, by emulating the Birman--Hilden theory for a certain class of manifolds called hyperk\"{a}hler manifolds.

\begin{definition}
A closed complex manifold $(M,J)$ of real dimension $4n$ is \textit{hyperk\"{a}hlerian} if there are two additional complex structures $I,K$ such that $\{I,J,K\}$ satisfy the quaternionic relations, and there exists a metric $g$ such that $(M,g)$ is K\"{a}hler with respect to $I$, $J$, and $K$. We say $M$ is \textit{hyperk\"{a}hler} if we assign to it the datum of such a compatible metric and complex structures. A hyperk\"{a}hler manifold $M$ is \textit{irreducible} if it is simply-connected with $\dim(H^{2,0}(M)) = \dim(H^0(M,\Omega^2))= 1$. 
\end{definition}
Unless otherwise stated, ``mapping class groups" will mean \textit{smooth} mapping class groups, so $\MOD(X) := \pi_0(\Diff^+(X))$ and $\SMod(X) := \SDiff_D^+(X)/\text{isotopy}$, where $\SDiff_D^+(X)$ is the normalizer of a subgroup $D \subset \Diff^+(X)$. The smooth mapping class group is one of the more mysterious variants of this object. For example, the following is a foundational problem in 4-manifold topology, as mentioned in \cite[pg. 12]{looijenga2021teichmuller}:

\begin{question}
Let $X$ be a K3 surface, and $\MOD(X) = \pi_0 (\Diff^+(X))$. Suppose $f \in \MOD(X)$ induces the identity action on $H^2(X,\ZZ)$. Is $f = \id \in \MOD(X)$? In other words, is the smooth Torelli group of $X$ trivial? 
\end{question}

Regarding \Cref{pi1Mod}, a connection between mapping class groups and moduli spaces in this hyperk\"{a}hler setting is mediated through \textit{Teichm\"{u}ller spaces}. This technology was famously leveraged by \cite{verbitsky2013mapping} in the proof of the global Torelli theorem for hyperk\"{a}hler manifolds. In response to Verbitsky's work, Kreck and Su were the first to give an example of a hyperk\"{a}hler manifold whose smooth Torelli group is infinite \cite{kreck2021finiteness}, or equivalently, one whose Teichm\"{u}ller space has infinitely many connected components. Recently, Looijenga gave a simplified proof of Verbitsky's global Torelli theorem \cite{looijenga2021teichmuller}. Pushing further on this Teichm\"{u}ller-theoretic perspective yields one of our main results:

\begin{mainthm}[Hyperk\"{a}hler Birman--Hilden Theorem]\label{IntroThmA}
Let $p: X \to Y$ be a smooth regular $D$-covering of a manifold $Y$ admitting a complex structure, where $X$ is a manifold admitting an irreducible hyperk\"{a}hler structure. Then the cover $p$ has the smooth Birman--Hilden property. That is, $\SMod(X)/D \cong \MOD(Y)$, where $\SMod(X)$ is the symmetric mapping class group and $\MOD(Y)$ is the mapping class group.
\end{mainthm}

The simplest example of a space covered by a hyperk\"{a}hler manifold is an Enriques surface:

\begin{definition}
Let $X$ be a K3 surface (a simply-connected compact complex surface which admits a nowhere vanishing holomorphic $2$-form). An \textit{Enriques surface} is a complex surface $Y$ biholomorphic to the quotient of $X$ by fixed-point free involution $\iota$.
\end{definition}

Here is a concrete special case of the main result that motivated this project. This Birman--Hilden type theorem allows us to pin down the smooth mapping class group of an Enriques surface up to isomorphism.

\begin{corollary}
Let $X$ be the smooth manifold underlying a K3 surface, $\iota : X \to X$ an involution without fixed points, and $Y = X/ \langle \iota \rangle$ its Enriques surface. Then $\SMod(X)/ \langle \iota \rangle \cong \MOD(Y)$, that is,  the cover associated to $\langle \iota\rangle$ has the smooth Birman--Hilden property.
\end{corollary}

In this respect, the present work is different from the articles of Hain \cite{hain2023mapping} and Kreck--Su \cite{kreck2021finiteness}, as our theorems hold even if $\dim_\RR(X) = 4$. The use of moduli spaces and their relationships to mapping class groups is much closer in spirit to the work of Farb--Looijenga on K3 surfaces (\cite{farb2021nielsen} and \cite{farb2023moduli}).

There is another distinct proof of the classical Birman--Hilden theorem due to MacLachlan--Harvey \cite{maclachlan1975mapping} via Teichm\"{u}ller theory. This is the model of proof we transfer to this setting, which forms the basis for half of this article. Let $Y$ to be a smooth complex manifold with hyperk\"{a}hlerian universal cover $X$ (such $Y$ are called \textit{Enriques manifolds}). The structure of our proof of the Birman--Hilden theorem, following \textit{loc. cit.}, is as follows: given the covering map $p :X \to Y$, there will be an associated lifting map between Teichm\"{u}ller spaces. Analyzing properties of this map, namely its injectivity, and related mapping class group actions will yield the main theorem. In fact, injectivity of this lifting map implies a global Torelli theorem for Enriques manifolds (see \Cref{TeichEnriques} for the definitions of the relevant spaces/maps involved):

\begin{mainthm}[Global Torelli Theorem for Enriques Manifolds]\label{IntroThmB}
Given an Enriques manifold $Y$, let $\Teich_s(Y)$ be its separated Teichm\"{u}ller space and $\mathcal{D}(Y)$ the period domain of complex structures on $Y$. The period mapping $\mathcal{P}_s: \Teich_s(Y) \to \mathcal{D}(Y)$ is injective on the connected components of $\Teich_s(Y)$.
\end{mainthm}

Part of our work will involve building the Teichm\"{u}ller spaces associated to Enriques manifolds, as we could not find this in the literature. Analogies between the moduli of K3 surfaces and Teichm\"{u}ller spaces of Riemann surfaces have appeared in many places, going at least as far back as Weil's ``Final Report" \cite{weil1958final}. To the best of my knowledge, this work appears to be the first attempt to construct Teichm\"{u}ller spaces of Enriques surfaces and their higher dimensional analogs, and to study various maps between Teichm\"{u}ller spaces in this setting. 

\subsection{Reviewing Birman--Hilden Theory}

We first recall the classical work of Birman--Hilden relating braid groups to the mapping class groups of hyperbolic surfaces, following \cite{FarbMargalit12}. Let $S_g^1$ be a surface of genus $g$ with one boundary component, and let $\MOD(S_g^1) = \pi_0(\Homeo^+(S_g^1))$ be its mapping class group. Let $\iota$ be a hyperelliptic involution of $S_g^1$, with $2g+1$ fixed points. Define \textit{the symmetric homeomorphism group} $\SHomeo^+(S_g)$ of $S_g^1$ to be the centralizer of $\iota$ in $\Homeo^+(S_g^1)$ (every map considered in this article will be orientation preserving).

The \textit{symmetric mapping class group} is the group
$$\SMod(S_g^1) = \SHomeo^+ (S_g^1)/ \text{isotopy},$$
that is, the subgroup of $\MOD(S_g^1)$ consisting of isotopy classes of homeomorphisms with symmetric representatives.

The quotient $S_g^1 / \langle \iota\rangle = D_{2g+1}$ yields a disk with $2g+1$ punctures, serving as the branch points of the cover. Every element of $\SHomeo^+(S_g^1)$ commutes with the involution $\iota$, and thus all such maps descend to homeomorphisms of the quotient disk. Moreover, commutativity also implies that such homeomorphisms preserve the set of $2g+1$ fixed points of $\iota$, and so there is a homomorphism
$$\SHomeo^+(S_g^1) \to \Homeo^+(D_{2g+1}).$$
This map is an isomorphism. 

Two symmetric homeomorphisms of a $S_g^1$ are \textit{symmetrically isotopic} if they are isotopic through symmetric homeomorphisms, i.e. they lie in the same path component of $\SHomeo(S_g^1)$. Given this definition, it follows that
\begin{align*}
    \SHomeo^+(S_g^1)/\text{symmetric isotopy} &= \pi_0(\SHomeo^+(S_g^1))\\
    &\cong \pi_0(\Homeo^+(D_{2g+1}))\\
    &= \MOD(D_{2g+1})\\
    &\cong B_{2g+1}.
\end{align*}    
The main goal is to show is that $\SMod(S_g^1)$ is isomorphic to the braid group $B_{2g+1}$. This amounts to showing that if two symmetric homeomorphisms of $S_g^1$ are isotopic, they can be made symmetrically isotopic with respect to the deck group $\langle \iota \rangle$. Covering maps that have this property are said to have the \textit{Birman--Hilden property}. Birman and Hilden proved this to be true \cite{birman1973isotopies}, strengthening the connection between the study of braid groups and mapping class groups of surfaces:

\begin{theorem}[Birman--Hilden Theorem]
Using the above notation, $\SMod(S_g^1) \cong B_{2g+1}$. There is analogous isomorphism $\SMod(S_g^2) \cong B_{2g+2}$, where $S_g^2$ is the genus $g$ surface with 2 boundary components.
\end{theorem}

\begin{remark}
The above isomorphism relates braid groups and their associated hyperelliptic mapping class groups by the isomorphism $B_{2g+1} \cong \SMod(S_g^1)$. Note that the symmetric mapping class group is not being quotiented by the deck group generated by the hyperelliptic involution, since the involution does not define an element of the mapping class group---$\iota$ does not fix the boundary pointwise. For the hyperk\"{a}hler Birman--Hilden theorem, quotienting $\SMod(X)$ by $D$ is necessary, as a nontrivial deck transformation represents a nontrivial symmetric mapping class.
\end{remark}

For a nice review of the classical Birman--Hilden Theorem and some generalizations, see the theses of Ghaswala \cite{ghaswala2017liftable} and Winarski \cite{winarski2015symmetry}; see also the survey article of Margalit--Winarski \cite{margalit2021braid}. Recent works of Ghaswala--Winarski \cite{ghaswala2017lifting} and McLeay \cite{mcleay2018subgroups} all provide new extensions of this result within the realm of surface topology, e.g. for totally ramified covers (covers where the preimage of a branch point is a point itself). There are still many interesting open questions left to answer, even in this low-dimensional setting.

\subsection{Why Hyperk\"{a}hler?}

Returning to the mapping class groups of $4$-manifolds, I was originally motivated to study the smooth mapping class group of Enriques surfaces $Y$. These are in some sense the simplest example of a non-simply connected complex surface, as their fundamental group is $\pi_1(Y) \cong \ZZ/2\ZZ$. 

In order to understand the mapping class group of an Enriques surface $Y$, it seems most natural to relate it to the mapping class group of its universal cover $X$, just as the Birman--Hilden theorem relates the mapping class group of a surface downstairs to its covering surface upstairs. Adapting a Birman--Hilden type theorem to this context in fact worked in the more general setting of hyperk\"{a}hler manifolds (which a K3 surface is an example of).

Hyperk\"{a}hler manifolds appear naturally in complex geometry, as they form one of the three fundamental building blocks of all K\"{a}hler manifolds $X$ with $c_1(X)_\RR = 0$. The Beauville-Bogomolov decomposition theorem says that all such $X$ admit a finite \'{e}tale covering by the product of complex tori, Calabi--Yau manifolds, and hyperk\"{a}hler manifolds (\cite{bogomolov1974decomposition}, \cite{beauville1983varietes}). Hyperk\"{a}hler geometry also arises when addressing fundamental questions in algebraic geometry, for example, questions regarding the rationality of cubic $4$-folds \cite{hassett1996special,beauville1985variete}. Quotients of hyperk\"{a}hler manifolds, called Enriques manifolds, were first introduced almost simultaneously in \cite{oguiso2011enriques} and \cite{boissiere2011higher}. They serve as the analog to Enriques surfaces in this higher dimensional setting. 

At the present time, hyperk\"{a}hler manifolds are some of the only known examples of higher dimensional spaces $X$ with a satisfactory Teichm\"{u}ller theory. Such marked moduli spaces are necessary in our analysis of higher dimensional mapping class groups. For example, ``the difference" between $\pi_1(\mathscr{M}_X)$ and $\MOD(X)$ is precisely recorded by the Teichm\"{u}ller spaces for hyperk\"{a}hler manifolds $X$; Torelli mapping classes can only arise in the monodromy groups of families varying $X$ over a \textit{disconnected} base $B$ \cite[Proposition 2.3]{looijenga2021teichmuller}. 

A stumbling block one runs into when trying to generalize the Birman--Hilden theory to higher dimensions is the lack of a result like the Dehn--Nielsen--Baer theorem, which informally tells us that a mapping class $f \in \MOD(S_g)$ is determined by its action on $\pi_1(S_g)$. The action of the Enriques mapping class group $\MOD(Y)$ on $\pi_1(Y)$ (a finite cyclic group, see \Cref{EnriquesReview}) reveals very little about $\MOD(Y)$. The next most obvious replacement, considering the $\MOD(Y)$-action on $H^2(Y,\ZZ)$, does not yield enough information to conclude the Birman--Hilden theorem, since the Torelli group completely evades analysis. Building the relevant Teichm\"{u}ller spaces using the ``Main Lemma" (see \Cref{TheMainLemma}) gives us just enough control to produce the desired isomorphism. 

\subsection{Nielsen Realization Results}
After proving the hyperk\"{a}hler analog of the Birman--Hilden theorem, we will come back down to the 4-manifold setting in order to resolve variants of the Nielsen realization problem for the Enriques surface. The key to establish these Nielsen realization results is the work of Farb--Looijenga \cite{farb2021nielsen}, which can be appealed to once lifted into the symmetric mapping class group of the K3 surface. Metric and complex Nielsen realization for K3 surfaces was settled in \textit{loc. cit.} via a homological criterion. This gives a sufficient and necessary linear-algebraic condition for a finite subgroup of the Enriques mapping class group to be metric/complex Nielsen realizable:

\begin{mainthm}[Enriques Metric/Complex Nielsen Realization]\label{IntroThmC}
Let $Y$ be an Enriques surface and $X$ its covering K3 surface. A finite subgroup $G \subset \MOD(Y)$ is realizable as a group of Ricci-flat isometries $G\subset \Isom(Y,h)$ (respectively K\"{a}hler-Einstein automorphisms $G \subset \Aut(Y,h,J)$) if and only if its group of lifts $\tilde{G} \subset \SMod(X)$ is Ricci-flat (respectively K\"{a}hler-Einstein) realizable on $X$.
\end{mainthm}

Through a similar line of reasoning, we are able to address smooth involutive realization for Enriques surfaces. For example, we prove smooth non-realizability of Dehn twists on Enriques surfaces:

\begin{mainthm}[Non-realizability of Dehn twists]\label{IntroThmD}
Let $T_S \in \Diff(Y)$ be the Dehn twist about an embedded $2$-sphere $S \subset Y$ with self intersection $S \cdot S = -2$. Although $T_S^2$ is smoothly isotopic to the identity, $T_S$ is not topologically isotopic to any finite order diffeomorphism of the Enriques surface $Y$.
\end{mainthm}

The techniques used in \cite{arabadji2023nielsen} give stronger non-realizability results on non-spin $4$-manifolds with finite fundamental group. In fact, their work already implies \Cref{IntroThmD}, and much more. However, this Birman--Hilden type theorem allows for a more elementary approach to proving smooth non-realizability.

One can ask how the metric and complex Nielsen realization problems compare on the Enriques surface. An example of a metric realizable subgroup of $\MOD(X)$ (isomorphic to the alternating group $A_4$) that preserves no complex structure on the K3 surface $X$ can be found in \cite[Theorem 1.4]{farb2021nielsen}. Thus the metric and complex Nielsen realization problems are distinct for K3 surfaces. The opposite holds true for Enriques surfaces $Y$, where we have the following:

\begin{mainthm}[Ricci-flat implies complex]\label{IntroThmE}
A finite subgroup $G \subset \MOD(Y)$ can be realized by isometries for some
Ricci-flat metric on the smooth Enriques surface $Y$ if and only if it preserves a compatible complex structure on $Y$.
\end{mainthm}

Any $G$-invariant complex structure $J$ on $Y$ can always be endowed with a compatible $G$-invariant Ricci-flat metric. To see this, suppose we are given any $J$-compatible K\"{a}hler class $\omega_0$ on $Y$. Averaging $\omega_0$ along the $G$-action yields a $G$-invariant K\"{a}hler class $\omega$, and we can then apply Yau's theorem to get the uniquely compatible $G$-invariant Ricci-flat metric associated to $\omega$. The novelty of \Cref{IntroThmE} is the ``only if" implication.

\subsection{Related Works}

The global Torelli theorem for Enriques manifolds (\Cref{IntroThmB}) crucially relies on the work of Verbitsky \cite{verbitsky2013mapping,verbitsky2020errata} and Looijenga \cite{looijenga2021teichmuller} in the hyperk\"{a}hler setting (see \cite{huybrechts2011global} for a survey of the former's work). In the case of K3 surfaces, Burns--Rapoport \cite{burns1975torelli}, Piatetski-Shapiro--Shafarevich \cite{pyatetskii1971torelli}, and Looijenga--Peters \cite{looijenga1980torelli} all gave (corrected) proofs of the global Torelli theorem using various algebraic/k\"{a}hlerian hypotheses before Siu proved all K3 surfaces are K\"{a}hler \cite{siu1983every}. The global Torelli theorem for Enriques surfaces was first shown by Horikawa \cite{horikawa1977periods} and refined by Namikawa \cite{namikawa1985periods}. Oguiso and Schr\"{o}er's papers on Enriques manifolds \cite{oguiso2011enriques} and their periods \cite{oguiso2011periods} lay the foundation for our proof of the global Torelli theorem in this higher dimensional setting.

As previously mentioned, our proof of the hyperk\"{a}hler Birman--Hilden theorem closely mirrors the classical argument of MacLachlan--Harvey \cite{maclachlan1975mapping}. Regarding the Nielsen realization results on Enriques surfaces, the work of Farb--Looijenga (\cite{farb2021nielsen} and \cite{farb2023moduli}) and Lee (\cite{lee2021nielsen} and \cite{lee2023isotopy}) have similarly resolved variants of the Nielsen realization problem for K3 surfaces and del Pezzo surfaces, respectively. Billi extended the work of Farb--Looijenga to the case of hyperk\"{a}hler manifolds \cite{billi2022note}. In \cite{baraglia2023note} and \cite{konno2022dehn}, Baraglia and Konno used Seiberg--Witten theory to show new instances of non-realizability which extend some of the results in \cite{farb2021nielsen} to general spin $4$-manifolds. Further non-realizability results for non-spin $4$-manifolds with finite fundamental groups recently appeared in the work of Arabadji--Baykur \cite{arabadji2023nielsen}, and for simply-connected indefinite non-spin $4$-manifolds by Konno--Miyazawa--Taniguchi \cite{konno2023involutions}. 

\subsection{Outline}
In \Cref{SmoothLifting}, we review some basic results on isotopy classes of diffeomorphisms and lifting them through covers. \Cref{TwistorReview} recalls the fundamentals of the twistor construction, as it will be an essential part of our construction of Teichm\"{u}ller spaces. \Cref{EnriquesReview} covers the basics of Enriques manifolds, discusses the known examples, and finishes with a lemma on extensions of Hodge isometries. In \Cref{PeriodDomains}, we introduce period domains and period mappings for Enriques manifolds. In \Cref{TeichEnriques}, we prove the ``Main Lemma" for Enriques manifolds, which allows us to build various Teichm\"{u}ller spaces and prove the Global Torelli theorem. \Cref{BHIsomorphism} proves the Birman--Hilden theorem for hyperk\"{a}hler manifolds using the lifting map between Teichm\"{u}ller spaces. In \Cref{NielsenSection}, we use the Birman--Hilden theorem to address various Nielsen realization problems for Enriques surfaces by lifting into the K3 setting. Finally, \Cref{TorelliConjecture} discusses a conjecture of Looijenga, which asserts that the smooth Torelli group of an Enriques surface is trivial. Our work from previous sections then shows that the smooth Enriques Torelli group is trivial if and only if the smooth K3 Torelli group is trivial (\Cref{iffTorelli}).

\subsection{Acknowledgements}
This work was supported in part by NSF Grant No. DMS-1944862. I am grateful to my advisor Jesse Wolfson, for his guidance and support throughout this project, and for giving extensive comments on many previous drafts. I would like to thank Benson Farb for suggesting the problem which eventually turned into this project, and for organizing a workshop on $4$-manifold mapping class groups which was very formative for me. Conversations with him and his comments on a previous write-up have greatly improved the exposition. I thank Ailsa Keating and Rebecca Winarski for their comments on an early draft of this paper. I thank David Baraglia, Tyrone Ghaswala, Richard Hain, Daniel Huybrechts, Hokuto Konno, Dan Margalit, and Curtis McMullen for their helpful comments and suggestions. I acknowledge Josh Jordan and Jeffrey Streets for many helpful conversations. I thank Eduard Looijenga for insightful discussions, and for his work which inspired many of these results. I thank the anonymous referee for numerous helpful comments and corrections.

\section{Smooth Lifting Lemmas}\label{SmoothLifting}

Much of this section is basic material that is well known, but since we shall often use these facts without reference, we include these statements for completeness. In what follows, our definition of the mapping class group will be taken in the smooth setting, as was mentioned in the introduction. In particular, all isotopies are smooth isotopies unless otherwise stated. The definitions are reviewed below. After setting up some notation, we discuss certain lifting lemmas that apply to all smooth manifolds. 

\subsection{Definitions and Lifting Isotopies}

\begin{definition}
For a given closed smooth manifold $M$, define its \textit{mapping class group} $\MOD(M)$ to be the group of smooth isotopy classes of orientation-preserving diffeomorphisms of $M$, i.e. 
$$\MOD(M) := \pi_0(\Diff^+(M)).$$
Suppose now there is a smooth group action by a deck group $D$ on $M$ that is proper and free. The symmetric diffeomorphism group $\SDiff_D^+(M)$ is the normalizer of $D$ in $\Diff^+(M)$. Define the \textit{symmetric mapping class group} of $M$ by
$$\SMod_D(M) := \SDiff_D^+(M)/\text{isotopy}.$$
When the deck group $D$ is clear from context, we suppress the subscript in the above notation, so then $\SDiff^+(M) = \SDiff_D^+(M)$ and $\SMod(M) = \SMod_D(M)$. Finally, consider the quotient manifold $N = M/D$. The \textit{liftable mapping class group} $\LMod(N)$ of $N$ is the subgroup of $\MOD(N)$ consisting of mapping classes that have representatives that lift to diffeomorphisms of $M$.
\end{definition}
Regular covers are endowed with the path lifting property, and so if we restrict our attention to the class of simply-connected manifolds $M$, all of $\MOD(N)$ lifts into $\pi_0(\Homeo(M))$, as every smooth mapping class of $N$ has a representative that lifts to a homeomorphism of $M$. We shall now check that each lift is smooth upstairs. 

\begin{proposition}
Let $M$ be simply-connected. The full mapping class group of $N$ will admit smooth representative lifts, i.e. $\MOD(N)=\LMod(N)$, through the cover $p: M \to N$ generated by $D$.
\end{proposition}
\begin{proof}
Let $\phi \in \MOD(N)$ and let $f$ be a diffeomorphism representative of $\phi$. We need only check that for every such diffeomorphism $f : N \to N$, its continuous lift $\tilde{f} : M \to M$ is smooth. Pick basepoints $e \in M$, $p(e) = x \in N$, and suppose $x = (f \circ p)(z)$. Recall that the generalized path lifting lemma gives us that
$$f \circ p : M \to N$$
lifts to a map $\widetilde{f\circ p} =: \tilde{f}:M \to M$ since $(f \circ p)_*(\pi_1(M,e)) \subset p_* (\pi_1(M,e))$, as $M$ is simply connected. Here is how we define the lift $\tilde{f}$. Given any $y\in M$, pick a smooth path $\gamma : [0,1] \to M$ with $\gamma(0) = z$ and $\gamma(1) = y$. Define 
$$\tilde{f}(y)= \text{the endpoint of the lift of } f \circ p\circ \gamma.$$
The fundamental group condition gives us that $\tilde{f}$ is well defined, as well as continuity of $f$ (note that we have connectedness and local path connectedness here, and that they are necessary). We can also quickly see that $\tilde{f}$ is a homeomorphism, for example by lifting $f \circ f\inv$.

To check smoothness of $\tilde{f}$, it suffices to locally represent $\tilde{f}$ as a composition of smooth maps. Let $u \in M$, then $(f\circ p)(u) \in N$. Since $p$ is a smooth covering map, there exists an open neighborhood $V \subset N$ containing $(f\circ p)(u)$ whose preimage $p\inv(V)$ is a disjoint collection of open subsets of $M$, each diffeomorphic to $V$. Consider the unique open subset $U \subset p\inv(V) \subset M$ characterized by the following features: $\tilde{f}(u) \in U$, and $U$ is diffeomorphic to $V$ via $p$. Then locally around $u \in \tilde{f}\inv(U)$ we have $\tilde{f} = p|_V \inv \circ (f\circ p)$, proving that $\tilde{f}$ is a smooth lift with smooth inverse. 
\end{proof}

\begin{lemma}\label{symlift}
Let $M$ be a simply connected closed manifold and let $p: M \to N$ a smooth regular cover with deck group $D$. Suppose $\Phi: N\times I \to N$ is a smooth isotopy between diffeomorphisms $f,h\in\Diff(N)$. Then $\Phi$ lifts to a smooth symmetric isotopy $\tilde{\Phi} : M \times I \to M$ between some lifted diffeomorphisms $\tilde{f}, \tilde{h} \in \SDiff(M)$ through $p$.
\end{lemma}
\begin{proof}
It suffices to prove the claim for when $h = \text{id}_N$, i.e. when $\Phi$ is a smooth isotopy between $f$ and the identity. 

Let $\Phi : N \times I \to N$ be so that $\Phi(x,0) = x$, $\Phi(x,1) = f(x)$, and $\Phi(\cdot,t)$ is a diffeomorphism for all other fixed $0 <t <1$. By generalized path lifting, simply-connectedness of $M\times I$ implies that there exists a lift of the isotopy
$$ \Tilde{\Phi}: M\times I \to M$$
between $\tilde{\id}_N$ and $\tilde{f}$, which locally looks like $\tilde{\Phi} = p \inv \circ (\Phi \circ (p \times \text{id}_I))$. Smoothness of $\tilde{\Phi}$ follows from the same logic as in the previous proof, and symmetry follow immediately from the path-lifting construction.
\end{proof}

\subsection{Reformulating the Birman--Hilden property}
We now want to study when a cover $p: M\to N$ has the Birman--Hilden property. First we begin with some generalities. Define a group homomorphism
$\Phi: \SDiff^+(M) \to \Diff^+(N)$ by
$$\Phi(\tilde{f})(x) = (p \circ \tilde{f})(y)$$
for any $y \in p\inv(x)$. This map is well-defined because all elements $\tilde{f} \in \SDiff^+(M)$ are fiber-preserving by definition. That is, whenever $p(a) = p(b)$ we have $(p\circ \tilde{f})(a) = (p \circ \tilde{f})(b)$. 

For any $\tilde{g} \in \ker(\Phi)$, $\Phi(\tilde{g}) (x) = x$ for all $x \in N$, and since $\tilde{g}$ must be fiber-preserving, this means that $\ker(\Phi) = D$. The image of $\Phi$ consists of all liftable diffeomorphisms $\LDiff^+(N)$, the group of diffeomorphisms on $N$ that have lifts who normalize the deck group action. Thus by the first isomorphism theorem, $\SDiff^+(M)/ D \cong \LDiff^+(N)$. 

We claim that having the Birman--Hilden property is equivalent to this isomorphism remaining true when we replace all the diffeomorphism groups involved with their respective isotopy groups, i.e. the associated mapping class groups. The following proposition makes precise this claim, and gives a group-theoretic characterization of coverings which enjoy the Birman--Hilden property. It is taken from the survey article \cite{margalit2021braid}. Note that as stated in \textit{loc. cit.} (without proof), this equivalence is specifically for quotients of hyperbolic surfaces $S$. Our proof of this more general statement is essentially identical in this higher dimensional case.

\begin{proposition}\label{BHequivalence}
Let $p: M\to N$ be a smooth regular covering map, generated by the finite deck group $D$. Furthermore, assume no two elements of $D$ are isotopic to one another. The following are equivalent:
\begin{enumerate}
    \item $p$ has the Birman–Hilden property: isotopic symmetric diffeomorphisms can be made symmetrically isotopic;
    \item The map $\LMod(N) \to \SMod(M)/D$ sending a liftable mapping class to its equivalence class of symmetric lifts is injective;
    \item The map $\Phi$ induces a well defined surjective map $\SMod(M) \to \LMod(N)$ sending a symmetric representative of symmetric mapping class downstairs;
    \item The map $\Phi$ induces an isomorphism $\SMod(M)/D \cong \LMod(N)$.
\end{enumerate}
\end{proposition}
\begin{proof}
($1 \Rightarrow 3$) Say the cover $p$ has the Birman--Hilden property. This is true if all $f$ and $g$ in $\SDiff^+(M)$ that are isotopic can be made symmetrically isotopic. This happens if and only if every equivalence class in $\SMod(M)$ can be defined up to symmetric isotopy. Thus there is an isomorphism 
$$\SDiff^+(M)/\text{symmetric isotopy} \cong \SMod(M).$$
The pushforward of a class $f \in \SMod(M)$ to $\LMod(N)$ via the map $\Phi$ thus gives rise to a well defined surjective group homomorphism. 

($3 \Rightarrow 4$) Since we assumed that no two elements of $D$ are isotopic, the kernel of the map $\SMod(M) \to \LMod(N)$ induced by $\Phi$ consists of the deck group $D$. Note that this is a slight abuse of notation, since $D$ also denotes the diffeomorphism subgroup, not just its isotopy classes. Thus we have the isomorphism $\SMod(M)/D \cong \LMod(N)$.

($4\Rightarrow 3$) This is immediately true, since the map $\SMod(M) \to \LMod(N)$ coming from $\Phi$ must be well-defined to induce the isomorphism $\SMod(M)/D \cong \LMod(N)$.

($2 \Leftrightarrow 4$) This is immediately true, as the map $\LMod(N) \to \SMod(M) / D$ is always surjective.

($3 \Rightarrow 1$) If the map $\SMod(M) \to \LMod(N)$ induced by $\Phi$ is well defined, that means that two representatives $f$ and $g$ of the same isotopy class in $\SMod(M)$ always land on the same element $[\Phi(f)] = [\Phi(g)] \in \LMod(N)$. Thus $\Phi(f)$ and $\Phi(g)$ are isotopic in $\LDiff(N)$. Any based lift of this isotopy through $p$ gives rise to a symmetric isotopy between $f$ and $g$ in $\SDiff(M)$, thus $p$ has the Birman--Hilden property.
\end{proof}

\section{A Review of The Twistor Construction}\label{TwistorReview}

We now review the basics of the twistor construction on irreducible hyperk\"{a}hler manifolds following \cite{looijenga2021teichmuller}, \cite{farb2023moduli}, and \cite{gross2012calabi}. It is central to building the relevant moduli spaces of geometric structures on the underlying smooth manifold $M$ of real dimension $4n$. Fix a Riemannian metric $g$ on $M$, let $H_R := H^2(M,\ZZ) \otimes_\ZZ R$ for any $\ZZ$-module $R$, and set $H = H_\ZZ$. The metric $g$ identifies the flat endomorphisms of the tangent bundle as a copy of the pure quaternions $\HH_{(M,g)}^{\mathrm{pure}} \subset \END(TM)$. Real Hodge theory then allows us to distinguish an oriented $3$-plane $P_g$ of harmonic $2$-forms in $H_\RR$ corresponding to an embedding of the flat endomorphisms $\HH_{(M,g)}^{\mathrm{pure}} \hookrightarrow H^2(M,\RR)$ \cite[pg. 3]{looijenga2021teichmuller}. 

Bogomolov, Beauville, and Fujiki showed there exists a nondegenerate, primitive, integral quadratic form $q: H \to \ZZ$ such that for some $c \in \QQ_{>0}$, the identity $q(a)^n = c \int_M a^{2n}$ holds for all $a \in H$ and for which $q_\RR$ is positive on the oriented $3$-plane $P_g$. The quadratic form $q$ is called the \textit{Beauville--Bogomolov--Fujiki form}. Given a holomorphic $2$-form $\omega$ on some hyperk\"{a}hler $X$ with $\int_X (\omega \bar{\omega})^n = 1$, the integral form $q(a)$ on $H_\RR$ can be expressed as a scalar multiple of
$$\frac{n}{2} \int_X a^2(\omega \bar{\omega})^{n-1} +(1-n)\left(\int_X a \omega^{n-1}\bar{\omega}^n \right) \left(\int_X a \omega^{n} \bar{\omega}^{n-1}\right).$$
Note that $q$ is in fact independent of any specific holomorphic $2$-form $\omega$ chosen. Furthermore, they proved the form $q_\RR$ has signature $(3,m)$, where $m = b_2(M) - 3$. This turns $H$ into an indefinite lattice of signature $(3,m)$; by a \textit{lattice}, we mean a free abelian group of finite rank endowed with an integral symmetric bilinear form. For the discussion that follows, the Grassmannian of positive $3$-planes in $H_\RR$ is denoted by $\Gr_3^+(H_\RR)$, and its associated symmetric metric $g_{loc}$. We assume that all of our manifolds $(M,g)$ have unit volume.

Given such a Riemannian metric $g$, the positive oriented $3$-plane $P_g \subset H_\RR$ it determines will orient the full tautological $3$-plane bundle over $\Gr_3^+(H_\RR)$, since $\Gr_3^+(H_\RR)$ is the (contractible) symmetric space associated to $\OO(q_\RR)$. This bundle orientation is called the \textit{spinor orientation}. We make the orientation and spinor orientation on $M$ part of our initial data.

If $J$ is a compatible complex structure on $M$ for which $g$ is K\"{a}hler, then the K\"{a}hler class $\kappa$ associated to $(M,g,J)$ lies in $P_g$ with 
$$\Vol_\kappa(M) := q_\RR(\kappa) =\int_M \sqrt{\det(g)} = 1.$$ 
The orthogonal complement of $\kappa$ in $P_g$ is a positive $2$-plane $\Pi$ that inherits a canonical orientation via the spinor orientation. This determines a complex line $H^{2,0}(M,J) \subset H^2(M,\CC)$: if $(x,y)$ forms an $\RR$-orthonormal basis for the $2$-plane $\Pi$, then the line $H^{2,0}(M,J)$ is spanned by $x +\sqrt{-1}y$ over $\CC$. Recall that $H^{2,0}(M,J)$ determines the Hodge decomposition on $H^{2}(M,\CC)$ uniquely.

If $g$ is also taken to be Ricci-flat then the converse holds as well, by Yau's theorem. Given any $\kappa \in P_g$ with $\Vol_\kappa(M) = 1$, there exists a complex structure $J$ on $M$ for which $g$ is a Ricci-flat K\"{a}hler metric and has $\kappa$ as its K\"{a}hler class. This complex structure is uniquely determined by the choice of $\kappa$, and all other such complex structures $J$ for which $g$ is a Ricci-flat K\"{a}hler metric also arise by choosing such a $\kappa \in P_g$ with $\Vol_\kappa(M)= 1$. We refer to the pair $(g,J)$ as a \textit{K\"{a}hler-Einstein structure} on $M$. Since $q_\RR$ is a positive quadratic form on the $3$-plane $P_g$, the space of complex structures compatible with a fixed unit volume Ricci-flat metric $g$ on $M$ is isomorphic to the 2-sphere 
$$\mathcal{D}(P_g) = \{\kappa \in P_g : \Vol_\kappa(M) = 1\} \cong S^2.$$

As we will build from scratch the Teichm\"{u}ller spaces attached to Enriques manifolds in \Cref{TeichEnriques}, we refer to \cite{looijenga2021teichmuller} for the definitions and constructions of Teichm\"{u}ller spaces and period mappings associated to hyperk\"{a}hler manifolds. A key difference between Looijenga's work versus others' found in the literature is the centrality of twistor families.

\section{Enriques Manifolds and The Kuranishi Family}\label{EnriquesReview}

\subsection{Cohomology and Families}

Let us now introduce the non-simply connected spaces to our story. Enriques manifolds were first introduced in \cite{oguiso2011enriques} and \cite{boissiere2011higher} with slightly different definitions. We will largely be following the treatment found in \cite{oguiso2011enriques} and \cite{oguiso2011periods}, especially regarding the construction of local period mappings.

\begin{definition}
    An \textit{Enriques manifold} is a connected complex manifold $Y$ that is not simply-connected, and whose universal covering $X$ is hyperk\"{a}hlerian.
\end{definition}

We shall now briefly describe all known examples of compact hyperk\"{a}hler manifolds, as they will be relevant in building the known examples of Enriques manifolds. Two infinite families of hyperk\"{a}hler manifolds, first introduced by Beauville, are constructed as follows:
\begin{enumerate}
    \item \textit{$\text{K3}^{[n]}$ manifolds}: the Hilbert scheme $S^{[n]}$ of $n$ points on a K3 surface $S$; this is a resolution of singularities of the quotient space $S^n / \mathfrak{S}_n$, and
    \item \textit{generalized Kummer manifolds}: given an abelian surface $A$, the Hilbert scheme of $n$ points $A^{[n+1]}$ is the minimal resolution of the symmetric product $A^{n+1}/\mathfrak{S}_{n+1}$. The latter space admits a ``summation map" to $A$, where $(z_0, \dots, z_n) \mapsto \sum z_i$, and so there is a composition map
    $$\varphi : A^{[n+1]} \to A^{n+1}/\mathfrak{S}_{n+1} \to A.$$
    The generalized Kummer manifold is then the fiber over $0\in A$, that is, $\mathrm{Kum}_n(A)= \varphi^{-1} (0).$
\end{enumerate}
There are two exceptional examples known, due to O'Grady (\cite{o1999desingularized} and \cite{o2000new}). In fact, all known compact hyperk\"{a}hler manifolds that have been constructed to date are deformation equivalent to one of the above.

We will often use $N$ to denote the smooth manifold underlying an Enriques manifold, and $\tilde{N}$ its smooth universal cover. We shall also use $Y$ to denote a complex Enriques manifold and $X$ its hyperk\"{a}hlerian universal cover. The trace of the representation $D = \pi_1(Y)$ on $H^{2,0}(X)$ gives rise to a homomorphism
$$D \to \CC^\times$$
which induces a canonical bijection $D \to \mu_d$, the multiplicative group of $d$-th roots of unity \cite[Lemma 2.3, Proposition 2.4]{oguiso2011enriques}. We thus will identify
$$D = \pi_1(Y) = \mu_d.$$
The integer $d \geq 2$ is called the \textit{index} of the Enriques manifold $Y$.

By \cite[Corollary 2.7]{oguiso2011enriques} or \cite[Proposition 6]{beauville1983some}, Enriques manifolds are necessarily projective, and so their universal covers are projective as well. Note that general hyperk\"{a}hler manifolds are not necessarily projective. 

The group of characters $\mu_d \to \CC^\times$ is cyclic of order $d$ with a canonical generator $\zeta \mapsto \zeta$. We identify
$$\Hom(\mu_d ,\CC^\times) = \ZZ/d\ZZ.$$
A finite-dimensional complex representation of $D$ is equivalent to a finite dimensional complex vector space $V$ endowed with a \textit{weight decomposition}
$$V = \bigoplus_{i \in \ZZ/d\ZZ} V_i$$
indexed by the character group. Explicitly, the \textit{weight spaces} $V_i \subset V$ is the subspace of vectors where each group element $\zeta \in D$ acts via multiplication by the complex number $\zeta^i \in \CC$. Note that $V_0 \subset V$ is the $D$-invariant subspace, and $V_1 \subset V$ is the subspace where the action of $\zeta$ is multiplication by itself.

Let $Y$ be an Enriques manifold of index $d\geq 2$, and $X \to Y$ its universal covering, so $X$ is a hyperk\"{a}hlerian manifold. The fundamental group $D$ acts on $H^1(X ,\Theta_X)$ such that we have a weight decomposition of cohomology vector spaces
$$H^1(X,\Theta_X) =\bigoplus_{i\in \ZZ/d\ZZ} H^1(X, \Theta_X)_i$$
indexed by the characters of $\pi_1(Y) = \mu_d$. Recall that $H^1(X,\Theta_X)_0$ is the $D$-invariant summand.

Let $\mathfrak{Y} \to B$ be a \textit{Kuranishi family} of $Y = \mathfrak{Y}_0$, that is, a deformation of $Y$ that is versal and has the property that $\dim(H^1(Y,\Theta_Y))=\dim(B)$.

\begin{proposition}[{\cite[Proposition 1.2]{oguiso2011periods}}]\label{EnriquesKuranishi}
    After shrinking $B$ if necessary, the Kuranishi family $\mathfrak{Y} \to B$ of an Enriques manifold $Y$ of index $d\geq 2$ is universal, the base is smooth, and each fiber $\mathfrak{Y}_b$ is an Enriques manifold of index $d$.
\end{proposition}
\begin{proof}[\proofname{ sketch}]
    There is a universal smooth Kuranishi family $\mathfrak{X}' \to B'$ for the hyperk\"{a}hler manifold $X = \mathfrak{X}_0$. By work of Fujiki \cite[Lemma 4.14]{fujiki1987rham} and Ran \cite[Corollary 2]{ran1992deformations}, after shrinking and modifying the family as need be to lie over a smooth base $B$, we can assume that $D = \pi_1(Y)$ acts on the family $\mathfrak{X} \to B$ fiberwise, which yields a family of index $d$ Enriques manifolds
    $$\mathfrak{X}/D \to B.$$
    Kodaira-Spencer theory and commutativity of the square
    \begin{center}
        \begin{tikzcd}
{H^1(Y,\Theta)} \arrow[r]       & {H^1(X,\Theta)}       \\
\Theta_B(0) \arrow[r] \arrow[u] & \Theta_B(0) \arrow[u]
\end{tikzcd}
    \end{center}
    implies versality of $\mathfrak{X}/D \to B$, and universality follows since $H^0(Y,\Theta_Y) = 0$ \cite[Proposition 1.1]{oguiso2011periods}.
\end{proof}

\subsection{Examples of Enriques Manifolds}

This section briefly describes the only known examples of Enriques manifolds $Y$ as first constructed in \cite{boissiere2011higher} and \cite{oguiso2011enriques}. We are most interested in understanding $H^2(Y,\ZZ)$ as a lattice, so that 
$$H^2(Y,\ZZ)_{\mathrm{free}} = H^2(Y,\ZZ)/\mathrm{torsion} \hookrightarrow H^2(X,\ZZ)$$
will naturally yield some primitive embedding of lattices. By transfer, the image of the injective map $H^2(Y,\ZZ)_{\mathrm{free}} \hookrightarrow H^2(X,\ZZ)$ is the $D$-invariant cohomology $H^2(X,\ZZ)^D$.

\begin{definition}
    Let $Y$ be an index $d$ Enriques manifold with universal $D$-cover $X$. A \textit{generalized Enriques lattice} $\Lambda_Y$ denotes the lattice given by 
    $$\Lambda_Y = \left(H^2(X,\ZZ)^D,q_r\right)$$
    where $q_r$ denotes the restriction of the quadratic form $q$ underlying the Beauville--Bogomolov--Fujiki lattice $\Lambda_X$ on $H^2(X,\ZZ)$.
\end{definition}

\begin{example}[From $\text{K3}^{[n]}$ manifolds]
    Consider a K3 surface $S$ with an Enriques involution $\iota$. The
Hilbert scheme $S^{[n]}$ comes with the natural involution $\iota^{[n]}$. If $n$ is odd, then $\iota^{[n]}$ has no fixed points, so
$$N = S^{[n]}/ \langle \iota^{[n]} \rangle$$
is an Enriques manifold of complex dimension $2n$ and index $2$. 

We shall now compute the lattice structure on the second cohomology of $N$. To do so, let us determine the invariant sublattice for the action of $\iota^{[n]}$ on $H^2(S^{[n]},\ZZ)$, which we equip with the Beauville--Bogomolov--Fujiki form. Recall that $H^2(S,\ZZ)^{\langle \iota \rangle} = U(2) \oplus E_8 (-2)$. It follows from naturality of the automorphism $\iota^{[n]}$ that
$$H^2(S^{[n]},\ZZ)^{\langle \iota^{[n]}\rangle} = U(2) \oplus E_8(-2) \oplus \langle -2(n-1) \rangle.$$
Let $p : S^{[n]} \to N$ denote the universal cover, $p^*$ the induced map on $H^2$, and $p_!$ the wrong-way map on $H^2$. By transfer, the composition $p_! \circ p^*: H^2(N,\ZZ)\to H^2(S^{[n]},\ZZ) \to  H^2(N,\ZZ)$ is multiplication by $2$. This computation, along with the universal coefficient theorem, implies that
$$H^2(N, \ZZ) = U \oplus E_8(-1) \oplus \langle -(n-1) \rangle \oplus \ZZ/2\ZZ.$$
There is another example one can build starting with a K3 surface using moduli of stable sheaves and Mukai vectors; we refer to \cite[Section 5]{oguiso2011enriques} for details, as the resulting cohomological computation yields an identical lattice on the second cohomology.
\end{example}

\begin{remark}
Note that the Hilbert scheme $(S/\langle \iota \rangle)^{[n]}$ does \textit{not} produce an Enriques manifold according to the definition found in \cite{oguiso2011enriques}, as its universal cover is not hyperk\"{a}hler. However, it is an Enriques manifold according to the definition found in \cite{boissiere2011higher}, as its universal cover is Calabi--Yau.
\end{remark}

\begin{example}[From generalized Kummer manifolds]
   Let $S$ be a bielliptic surface, that is, a surface $S$ with torsion canonical class of order $d \in \{2, 3, 4\}$ (here we omit the case of $d=6$ since it does not yield an Enriques manifold). It admits a finite \'{e}tale covering by an abelian surface $A \to S$. Certain conditions on $d$ and $n$, most notably that $d | (n+1)$, ensure that $D = \ZZ/d\ZZ$ acts on $A^{[n+1]}$ in such a way that it leaves $\mathrm{Kum}_n(A) \subset A^{[n+1]}$ invariant. In fact, this $D$-action on $\mathrm{Kum}_n(A)$ will be free. The quotient manifold $K = \mathrm{Kum}_n(A)/D$ is then an Enriques manifold of index $d$. Let $g_d \in D$ be a generator.
   
   One can show that the lattice structure on $H^2(\mathrm{Kum}_n(A),\ZZ)$ equipped with the Beauville--Bogomolov--Fujiki form is given by
   $$H^2(\mathrm{Kum}_n(A),\ZZ) =  U \oplus U \oplus U\oplus \langle -2(n+1)\rangle.$$
   where $H^2(A,\ZZ) = U \oplus U \oplus U$. We would like to determine the invariant sublattice of $H^2(\mathrm{Kum}_n(A),\ZZ)$ under the $D$-action. By naturality of the $D$-action on $\mathrm{Kum}_n(A)$, the summand $\langle -2 (n+1)\rangle$ lies in the invariant sublattice. To determine the remaining invariant pieces, we need to study the $D$-action on $H^2(A,\ZZ)$. According to the table found in \cite[pg. 15]{oguiso2011enriques}, it suffices to take $A = E \times F$, where $E$ and $F$ are specific CM elliptic curves classified by Bagnera and de Franchis. We briefly recall the classification now.
   
   Let $E = \CC/ (\ZZ + \tau_1 \ZZ)$ and $F = \CC/(\ZZ+ \tau_2\ZZ)$, where $\tau_1,\tau_2 \in \HH^2$ are the periods of $E$ and $F$, respectively. Let $\zeta = \exp(2\pi i/3)$ and $z \in F$ be arbitrary. Let $\tau_1$ be arbitrary, and set $\tau_2 = -1\text{, }\zeta, \text{ or } i$ when $d = 2\text{, }3 \text{ or }4$, respectively. Consider the following standard $D$-actions on $A = E \times F$ written in the coordinates $(e,f)$:
   \begin{align*}
       \ZZ/ 2\ZZ : \quad g_2: (e,f) &\longmapsto (e + 1/2, -f + z),\\
       \ZZ/ 3\ZZ : \quad g_3: (e,f) &\longmapsto (e + 1/3 , \zeta f + z),\\
       \ZZ/ 4\ZZ : \quad g_4: (e,f) &\longmapsto (e + 1/4 ,  i f + z).
   \end{align*}
   Since translation is cohomologically trivial, we need only to consider the multiplication part to determine the induced automorphism of $H^2(A,\ZZ)$. On $E$, let us write the coordinate $z = z_1 + i z_2$, and on $F$ we write $w = w_1+ iw_2$. Then $H^2(A,\RR)$ is generated by the $2$-forms 
   $$dz_1\wedge dw_1, dz_1\wedge dw_2,dz_2\wedge dw_1,dz_2\wedge dw_2,dz_1\wedge dz_2,dw_1\wedge dw_2.$$
   We can then explicitly compute the invariant sublattice in each case. We see that
   \begin{align*}
       H^2(\mathrm{Kum}_n(A),\ZZ)^{\langle g_2 \rangle} = U(2)\oplus \langle -2(n+1) \rangle,\\
       H^2(\mathrm{Kum}_n(A),\ZZ)^{\langle g_3 \rangle} = U(3)\oplus \langle -2(n+1) \rangle,\\
       H^2(\mathrm{Kum}_n(A),\ZZ)^{\langle g_4 \rangle} = U(4)\oplus \langle -2(n+1) \rangle.
   \end{align*}
   It follows from transfer and the universal coefficient theorem that in each case,
   $$H^2(K,\ZZ) = U \oplus \langle -2(n+1)/d \rangle \oplus \ZZ/d\ZZ.$$
   We refer the interested reader to \cite[Section 6]{oguiso2011enriques} and \cite[Section 4.2]{boissiere2011higher} for more details on building this quotient $K = \mathrm{Kum}_n(A)/D$.
\end{example}

Having studied all known examples, we now finish this section with the following lemma inspired by \cite[Theorem 1.4]{namikawa1985periods}, which will be useful in characterizing extensions of Hodge isometries:

\begin{lemma}\label{HodgeExtends}
    Let $\Lambda_Y$ denote one of the known generalized Enriques lattices for $Y$ an Enriques manifold with universal $D$-cover $X$, and $\Lambda_{X}$ the Beauville--Bogomolov--Fujiki lattice of $X$. Given any two primitive lattice embeddings $j_1,j_2 : \Lambda_{Y} \to \Lambda_{X}$, any isometry $\phi : \Lambda_{Y} \to \Lambda_{Y}$ extends to an isometry $\tilde{\phi} : \Lambda_{X} \to \Lambda_{X}$ satisfying $\tilde{\phi} \circ j_1 = j_2 \circ \phi$.
\end{lemma}
\begin{proof}
    We need only concern ourselves with the even unimodular part of the lattice, as any primitive embedding into the sublattice coming from the exceptional divisor, i.e. the sublattice $\langle -2 (n\pm 1)\rangle$, is an isomorphism. Regarding the even unimodular parts, these follow exactly the exact same argument outlined in \cite{namikawa1985periods}; the fundamental results of Nikulin (\cite[Propositions 1.5.1 and 1.14.1]{nikulin1980integral}) apply here and so we are done.
\end{proof}

\begin{remark}
    Note that an isometry of Hodge structures does not always imply that the underlying manifolds are biholomorphic or bimeromorphic to one another; in fact, this is false for general hyperk\"{a}hler manifolds.
\end{remark}

\section{Period Domains and Local Torelli}\label{PeriodDomains}

Our next task is to define period domains and period maps for Enriques manifolds, again following \cite{oguiso2011periods}. We shall indicate where our treatment diverges from theirs. Let $Y$ be a complex Enriques manifold and $X$ its hyperk\"{a}hlerian universal cover. Since $X$ has been given a complex structure, we can select a nondegenerate holomorphic $2$-form $\omega$ on $X$. The second cohomology $H = H^2(X,\ZZ)$ can then be equipped with the primitive integral Beauville--Bogomolov--Fujiki form $q_\omega$, for which $D$ acts on $(H,q_\omega)$, producing a $q_\omega$-orthogonal representation.

We will say a lattice has signature $(p,*)$ if it is signature $(p,q)$ for some $q \geq 0$.

When convenient, we shall let $q(a,b) = (a,b)$ be the bilinear extension of $q$ to $H^2(X,\CC)$ and $q_\CC (a,b) = (a,\Bar{b})$ its Hermitian extension.

\begin{proposition}[{\cite[Lemma 2.1]{oguiso2011periods}}]\label{latticeSignature}
    The lattice $(H^2(X,\ZZ),q)$ is nondegenerate of signature $(3,*)$. The Hermitian extension $q_\CC$ on the weight space $H^2(X,\CC)_1$ is nondegenerate, and has signature $(2,*)$ for $d = 2$, and $(1,*)$ for $d \geq 3$.
\end{proposition}

We shall now diverge from the work found in \cite{oguiso2011periods} for our definition of markings:

\begin{definition}
     A \textit{marking} of the complex Enriques manifold $Y$ is the choice of a $D$-invariant diffeomorphism $\Tilde{\varphi} : \tilde{N} \to X$, where $N$ is the smooth manifold underlying $Y$, and $X$ is the hyperk\"{a}hlerian manifold $D$-equivariantly diffeomorphic to the universal cover $\tilde{N}$. Pushing forward complex structures through the universal coverings $\tilde{N} \to N$ and $X \to Y$ yields a geometric marking downstairs $\varphi: N \to Y$.
\end{definition}

A $D$-invariant marking $\tilde{\varphi}$ of the manifold $X$ (as defined above), induces a \textit{lattice marking} of the second cohomology group $\tilde{\varphi}^* \! : \! (H^2(X,\ZZ),q_{\omega}) \to (H^2(\tilde{N},\ZZ), q)$ endowed with an orthogonal representation of $D$. Moreover, the weight space $H^2(\tilde{N},\CC)_1 = H_{\CC,1}$ is $q_\CC$-nondegenerate of signature $(2,*)$ for $d = 2$ and $(1,*)$ for $d \geq 3$.

\begin{remark}
    Looijenga points out that moduli spaces of lattice marked hyperk\"{a}hlerian manifolds seem to be preferred in the literature \cite[pg. 12]{looijenga2021teichmuller}. We choose to work with geometric markings here as they are essential in the construction of Teichm\"{u}ller spaces, which are central to the work of Verbitsky \cite{verbitsky2013mapping}, Markman \cite{markman2021existence}, and Looijenga \cite{looijenga2021teichmuller} on the Torelli theorems and construction of universal families for hyperk\"{a}hler manifolds. This is because geometric markings carry topologically information that lattice markings neglect, e.g. the spinor orientation. 
\end{remark}

\begin{definition}
    The \textit{period domain} $\mathcal{D}(Y)$ of marked Enriques manifolds $Y$ is
    $$\mathcal{D}(Y) := \{[\omega] \in \PP(H_{\CC,1}) : (\omega,\omega) = 0 \quad\text{ and }\quad (\omega,\overline{\omega} )>0)\},$$
    where $\overline{\omega}$ denotes complex conjugation inside the complexification $H_{\CC}$.
\end{definition}

Note that for $d = 2$, the weight space $H_{\CC,1} \subset H_{\CC}$ is invariant under complex conjugation. When $d \geq 3$, each $\omega \in H_{\CC,1}$ satisfies
$$(\omega,\omega) = (\zeta \omega,\zeta\omega) = \zeta^2 (\omega,\omega)$$
for all $\zeta \in G$, hence the weight space $H_{\CC,1} \subset H_{\CC}$ is totally isotropic; now the period domain is given by
$$\mathcal{D}(Y) := \{[\omega] \in \PP(H_{\CC,1}) :(\omega,\overline{\omega} )>0)\}.$$
Since our period domains are open subsets of a smooth quadric in projective space, with respect to the classical topology, they inherit the structure of a complex manifold.

\begin{proposition}[{\cite[Proposition 2.2]{oguiso2011periods}}]
    Set $k = \dim(H_{\CC,1})$. For $d = 2$, the period domain $\mathcal{D}(Y)$ is the disjoint union of two copies of bounded symmetric domains of Type $\mathrm{IV}_{k-1}$ of dimension $k-1$.

    For $d \geq 3$, the period domains $\mathcal{D}(Y)$ are bounded symmetric domains of type $\mathrm{I}_{1,k-1}$, and thus are biholomorphic to the complex ball of dimension $k-1$.
\end{proposition}

Given a marked complex Enriques manifold $(Y, \tilde{\varphi}\! : \!\tilde{N} \to X)$, the unique up to scaling holomorphic $2$-form $\omega \in H^{2,0}(X)$ satisfies
$$(\omega ,\omega) = 0 \quad \text{and} \quad (\omega,\overline{\omega}) >0,$$
with $\omega \in H^2(X,\CC)_1 \subset (H^2(X,\CC), q_\omega)$. Thus we can define the \textit{period point} of our marked Enriques manifold $(Y,\tilde{\varphi})$ to be $[\tilde{\varphi}^* (\omega)] \in \mathcal{D}(Y)$.

Let $\pi: \mathfrak{Y} \to B$ be a flat family of complex Enriques manifolds over a simply-connected complex base $B$. By universality of the Kuranishi family, \Cref{EnriquesKuranishi} tells us that each fiber $\mathfrak{Y}_b$ is an Enriques manifold of index $d$. Moreover, the universal covering $\mathfrak{X} \to \mathfrak{Y}$ is fiberwise the universal covering, and so we obtain a flat family $\tilde{\pi}:\mathfrak{X} \to B$ of hyperk\"{a}hlerian manifolds.

Suppose we are given a marking $\Tilde{\pi}_0 : \tilde{N} \to X$, where $X = \mathfrak{X}_0$ is the universal cover of $Y = \mathfrak{Y}_0$. Since the local system $R^2 \tilde{\pi}_* \ZZ_{\mathfrak{X}}$ is constant, simply-connectedness of $B$ gives us that the induced lattice marking of $Y$ given by $\tilde{\pi}_0^*$  \textit{uniquely} extends to a $D$-invariant lattice marking $\Tilde{\pi}^* : R^2 \tilde{\pi}_* \ZZ_{\mathfrak{X}} \to H_B$ over the flat family of Enriques manifolds. 
In turn, we obtain a period mapping
\begin{align*}
    \mathcal{P} : B &\longrightarrow \mathcal{D}(Y)\\
    b &\longmapsto [\tilde{\pi}^*(\omega_{\mathcal{X}_b})]
\end{align*}

By work of Griffiths \cite{griffiths1968periods}, the map $\mathcal{P}$ is holomorphic. We can now state the local Torelli theorem:

\begin{theorem}[Local Torelli Theorem, {\cite[Theorem 2.4]{oguiso2011periods}}]\label{localEnriques}
    Let $Y$ be a marked Enriques manifold and $\mathfrak{Y} \to B$ the Kuranishi family of $Y = \mathfrak{Y}_0$. Then the period map $\mathcal{P} : B \to \mathcal{D}(Y)$ is a local isomorphism.
\end{theorem}

\begin{remark}
    One can use the local Torelli theorem as in \cite[Definition 25.4]{gross2012calabi} to construct a course moduli space of isomorphism classes of marked Enriques manifolds $\mathscr{M}_Y$, which in turn yields a global period map into $\mathcal{D}(Y)$ which is \'{e}tale. However, one should be careful, as $\mathscr{M}_Y$ is not Hausdorff; \cite[Theorem 4.1]{oguiso2011periods} tells us that birational Enriques manifolds will have the same period point in $\mathcal{D}(Y)$, just as with hyperk\"{a}hler manifolds. The Teichm\"{u}ller space of Enriques manifolds that we will construct is also not separated.
\end{remark}

\section{Teichm\"{u}ller Spaces of Enriques Manifolds}\label{TeichEnriques}

\subsection{The Main Lemma}
We shall now approach the construction of the Teichm\"{u}ller space of Enriques manifolds, starting with the following lemma on families of Enriques manifolds. The statement of the lemma and pieces of the proof are modelled nearly verbatim off of \cite[Proposition 2.3]{looijenga2021teichmuller}. Looijenga remarks that the archetypical version of this result is the ``Main Lemma" of Burns--Rapoport \cite{burns1975torelli}. 

\begin{lemma}\label{TheMainLemma}
Let $\pi : \mathfrak{Y} \to U$ and $\pi': \mathfrak{Y}' \to U$ be proper holomorphic families of Enriques manifolds $Y$ and $Y'$ over the same simply-connected complex manifold $U$. Suppose we are given an isomorphism between the associated integral variations of Hodge structure (modulo torsion) in degree two:
$$\phi : R^2 \pi_*' \ZZ_{\mathfrak{Y}'} \cong R^2 \pi_* \ZZ_{\mathfrak{Y}}.$$
If for some $p \in U$, $\phi_p$ is induced by an isomorphism $f_p : Y_p \cong Y_p'$, then there exists a proper generically finite morphism $\hat{U} \to U$, a closed analytic subspace $\mathfrak{Z} \subset \mathfrak{Y} \times_{\hat{U}} \mathfrak{Y}'$ flat over $\hat{U}$, and a closed analytic, proper subset $K \subsetneq \hat{U}$ such that

\begin{enumerate}
    \item if $\hat{u} \in \hat{U}\backslash K$ lies over $u \in U$, then $Z_{\hat{u}}$ is the graph of an isomorphism $f_{\hat{u}} : Y_u \cong Y_u'$ which is isotopic to $f_p$ and induces $\phi_u$. Moreover, $f_p$ appears in this manner: for some $\hat{p} \in \hat{U} \backslash K$, $f_{\hat{p}} = f_{p}$,
    
    \item for every $u \in U$, $Y_u$ and $Y_u'$ are bimeromorphically equivalent,
    
    \item if there exists $\kappa \in H^0(U, R^2 \pi_* \RR)$ and $\kappa ' \in H^0(U, R^2 \pi_*' \RR)$ which restrict to a K\"{a}hler class in every fiber of $\pi$, respectively $\pi'$, and $\phi_p(\kappa'(p)) = \kappa (p)$, we can take $\hat{U} = U$ and $\mathfrak{Z}$ will be the graph of a $U$-isomorphism $\mathfrak{Y} \cong \mathfrak{Y}'$,
    
    \item the group $\Aut_0(\mathfrak{Y}/U)$ of automorphisms of $\mathfrak{Y}/U$ that are fiberwise isotopic to the identity is finite. Furthermore, it specializes for every $u\in U$ to the group $\Aut_0(Y_u)$ of automorphisms of $Y_u$ isotopic to the identity, and is via $f_p$ naturally identified with $\Aut_0(\mathfrak{Y}'/U)$. 
\end{enumerate}    
\end{lemma}
\begin{proof}
    We are given the proper holomorphic families of Enriques manifolds $\pi : \mathfrak{Y} \to U$ and $\pi': \mathfrak{Y}' \to U$ of $Y = \mathfrak{Y}_0$ and $Y' = \mathfrak{Y}_0'$, respectively. The proof of \Cref{EnriquesKuranishi} tells us that taking the universal covers of each family $\mathfrak{X} \to \mathfrak{Y}$ and $\mathfrak{X}' \to \mathfrak{Y}'$ yields $D$-invariant proper holomorphic families $\Tilde{\pi} : \mathfrak{X} \to U$ and $\Tilde{\pi}' : \mathfrak{X}' \to U$ of hyperk\"{a}hlerian manifolds $X$ and $X'$ that fiberwise cover $Y$ and $Y'$.
    
    We would like the isomorphism $\phi$ between integral variations of Hodge structure associated to the Enriques families to extend to a $D$-invariant isomorphism of the local systems of Beauville--Bogomolov--Fujiki lattices
    $$\Tilde{\phi} : R^2 \pi_*' \ZZ_{\mathfrak{X}'} \cong R^2 \pi_* \ZZ_{\mathfrak{X}}.$$
    By assumption, the isometry $\phi_p$ was induced by the isomorphism $f_p : Y_p \to Y_p'$, which then lifts to a $D$-invariant isomorphism of $\tilde{f}_p : X_p \to X_p'$. This then induces a $D$-invariant isometry of Beauville--Bogomolov--Fujiki lattices 
    $$\tilde{f}_p^* : H^2(X_p',\ZZ)\to H^2(X_p,\ZZ)$$ which preserves K\"{a}hler cones and holomorphic $2$-forms, so it is a Hodge isometry which by construction agrees with the $D$-invariant lattice isometry $\tilde{\phi}_p$. Simply-connectedness of the base $U$ allows us to conclude there exists some $\tilde{\phi}$ that is a $D$-invariant isomorphism of the integral variations of Hodge structures associated to $\mathfrak{X}/U$ and $\mathfrak{X}'/U$, by extending $\tilde{f}_p^*$. By construction, this $D$-invariant isomorphism of local systems $\tilde{\phi}$ induces $\phi$. 
    
    By \cite[Proposition 2.3]{looijenga2021teichmuller}, properties (1), (2), (3), and (4) hold for the families $\mathfrak{X}/U$ and $\mathfrak{X}'/U$, but to prove our lemma requires showing that this descends to the underlying families of Enriques manifolds. This amounts to checking that $D$-invariance of the families $\mathfrak{X}/U$ and $\mathfrak{X}'/U$ behaves nicely with each property listed. Note that the fiber product $\mathfrak{X} \times_U \mathfrak{X}'$ lies over $U$ and is equipped with the natural diagonal $D$-action.

    It will be helpful to recall parts of the proof of the Main Lemma for hyperk\"{a}hlerian families. Let $B = \mathscr{D}_{\mathfrak{X} \times_U \mathfrak{X}'}$ be the relative Douady space which parametrizes compact analytic subspaces of $\mathfrak{X} \times_U \mathfrak{X}'$ contained in a fiber of $\mathfrak{X} \times_U \mathfrak{X}'/U$. By work of Pourcin \cite{pourcin1969theoreme}, this is an analytic space, and comes with a universal family $\tilde{\mathfrak{Z}}_B \subset \mathfrak{X} \times_U \mathfrak{X}' \times_U B$ that is proper and flat over $B$. Let $\hat{U}$ be the irreducible component of $B$ which contains the $D$-invariant graph of $\tilde{f}_p : X_p \cong X_p'$ and set $\tilde{\mathfrak{Z}} = \tilde{\mathfrak{Z}}_{\hat{U}}$. The projection map of the family will be denoted $\hat{\pi} : \tilde{\mathfrak{Z}} \to \hat{U}$. Set-theoretically, we have
    $$\tilde{\mathfrak{Z}} = \{(X_{\hat{u}},X_{\hat{u}}') \in \mathfrak{X} \times \mathfrak{X}' : \hat{\pi}(X_{\hat{u}}) = \hat{\pi}(X_{\hat{u}}') =\hat{u}\in \hat{U}\}.$$
    Since $\tilde{\mathfrak{Z}}$ inherits the diagonal fiberwise $D$-action which fixes the base $\hat{U}$, it descends to a family $\mathfrak{Z}$ of Enriques manifolds. By studying $D$-invariant properties of this family, we shall prove the claim.
    
    The map $\mathfrak{X} \times_U \mathfrak{X}' \to U$ is a \textit{weakly K\"{a}hler} morphism, meaning there is a $2$-form on the source whose restriction to every fiber is a K\"{a}hler form. By work of Fujiki \cite{fujiki1987rham}, we then know that the projection $r : \hat{U} \to U$ is proper, and so $r(\hat{U})$ is a closed subvariety of $U$. We will show that it is all of $U$.
    
    By the local Torelli theorem for Enriques manifolds (\Cref{localEnriques}), there is an open neighborhood $V$ of $p$ in $U$ such that $f_p$ extends to a $V$-isomorphism $\mathfrak{Y}_V \to \mathfrak{Y}_V'$. This lifts to a $D$-invariant $V$-isomorphism of the covering families $\mathfrak{X}_V \to \mathfrak{X}_V'$. By construction, the graph of this isomorphism appears in $\tilde{\mathfrak{Z}} = \tilde{\mathfrak{Z}}_{\hat{U}}$. Because $V$ lies in the image of $r$, Zariski density of $V$ in $r(\hat{U})$ lets us conclude that $r(\hat{U}) = U$. Moreover, the locus $K$ of $ \hat{u}\in \hat{U}$ for which $\tilde{Z}_{\hat{u}}$ is not the graph of an isomorphism is a proper closed analytic subset of $\hat{U}$, and $r(K)$ is a proper closed analytic subset of $U$.

    The above shows that any $\hat{u} \in \hat{U} \backslash K$ lying over $u \in U$ has $\tilde{Z}_{\hat{u}}$ equal to the graph of a $D$-invariant isomorphism $\tilde{f}_{\hat{u}} : X_{u} \cong X_{u}'$ which is isotopic to $\tilde{f}_p : X_p \cong X_p'$ and induces $\tilde{\phi}_u$. To push this downstairs onto the families $\mathfrak{Y}$ and $\mathfrak{Y}'$, we need to show that this isotopy $\tilde{\Phi}$ can be promoted to a $D$-equivariant isotopy. This is clear given the construction of the universal family $\tilde{\mathfrak{Z}}$: it is a pullback of families of $D$-invariant hyperk\"{a}hlerian manifolds lying over the irreducible base $\hat{U}$. Moreover, since $K \subset \hat{U}$ is Zariski closed, its complement $\hat{U}\backslash K$ is path connected in the classical topology. An isotopy between the $D$-invariant isomorphism $\tilde{f}_u$ and $\tilde{f}_p$ can be represented by a path $\gamma$ in the base $\hat{U}\backslash K$. Geometrically we are deforming the graphs of such isomorphisms $\tilde{Z}_{u}$ and $\tilde{Z}_{p}$ to one another along the path $\gamma$. Any fiber lying over a point $q \in \gamma$ is the graph of another $D$-invariant isomorphism $\tilde{Z}_q$. Thus the isotopy between $\tilde{f}_{\hat{u}}$ and $\tilde{f}_p$ represented by $\gamma$ is taken through $D$-invariant graphs $\tilde{Z}_q$, proving that $\tilde{\Phi}$ is $D$-equivariant, and so $\Tilde{\Phi}$ descends to an isotopy $\Phi$ between $f_{\hat{u}}$ and $f_{p}$. This gives us Property (1). For Property (2), the exact same (now standard) proof given in \cite{looijenga2021teichmuller} works. 

    Let us show Property (3) now. By assumption, there exists global $2$-forms which restrict fiberwise to K\"{a}hler classes, and that the Hodge isometry $\phi_p$ must send a K\"{a}hler class to a K\"{a}hler class. Then the $\mathfrak{Y}/U$ and $\mathfrak{Y}'/U$ are families of K\"{a}hler Enriques manifolds, so their $D$-invariant covering families $\mathfrak{X}/U$ and $\mathfrak{X}'/U$ are also such families of hyperk\"{a}hlerian manifolds with $\tilde{\phi}_p$ sending a K\"{a}hler class to a K\"{a}hler class. Therefore the assumptions of Property (3) in \cite[Proposition 2.3]{looijenga2021teichmuller} are satisfied. We can conclude then that $K = \emptyset$ and $\hat{U} = U$.

    Finally, the proof of Property (4) follows the exact same line of reasoning as in the hyperk\"{a}hler setting: the (finite) group of isotopically trivial automorphisms of $Y$ fixes a K\"{a}hler-Einstein metric, and we apply Property (3) to two copies of $\mathfrak{Y}/U$ to get an extension to the full family.
\end{proof}

\subsection{Teichm\"{u}ller Spaces and Global Torelli Theorems}
Given \Cref{TheMainLemma}, we can now define and topologize a Teichm\"{u}ller space of Enriques manifolds.

\begin{definition}
    The \textit{Teichm\"{u}ller space} of an Enriques manifold $N$ is
$$\Teich(N) = \{\text{marked complex structures on }N\}/\Diff_0(N),$$ 
i.e. the moduli space of marked smooth isotopy classes of complex structures $[Y, \tilde{\varphi} : \tilde{N} \to X]$.
\end{definition}

As of right now, all we have defined is a Teichm\"{u}ller \textit{set}. Let us build an atlas of charts on $\Teich(N)$, following \cite[pg. 9-10]{looijenga2021teichmuller}:

Given an open subset $U$ of $\mathcal{D}(Y)$, a basic chart for $\Teich(N)$ with domain $U$ is given by a complex structure on
$N\times U$ for which the resulting complex manifold $\mathfrak{Y}$ has the property that
\begin{enumerate}
    \item the projection $\mathfrak{Y} \to U$ is holomorphic,
    \item the fibers of $\mathfrak{Y} \to U$ are complex Enriques manifolds, and
    \item its period map is given by the inclusion of $U$ into $\mathcal{D}(Y)$.
\end{enumerate}
Such an object will define an injection $U\hookrightarrow \Teich(N)$. The local Torelli theorem (\Cref{localEnriques}) tells us that every complex structure on $N$ appears as a member of such a family, so the basic charts cover all of $\Teich(N)$. We give $\Teich(N)$ the quotient topology, that is, the coarsest topology for which all the basic charts are continuous. It follows from \Cref{TheMainLemma} (with $\phi$ the identity and $f_p$ isotopic to the identity) that the locus where two basic charts with domains $U$ and $U'$ of $\mathcal{D}(Y)$ agree is the complement of a closed (analytic) subset of $U\cap U'$. Thus every basic chart is an open map. This makes our atlas complex-analytic and gives $\Teich(N)$ the structure of a (non-separated) complex manifold for which $\mathcal{P}$ is a local isomorphism.

It will be useful to define related Teichm\"{u}ller spaces, as Looijenga does \cite{looijenga2021teichmuller}, to state and prove analogous refined versions of the global Torelli theorem. To start off, \cite[Theorem 4.1]{oguiso2011periods} tells us that, just as in the hyperk\"{a}hler case, $\Teich(N)$ is not Hausdorff. This motivates the following definition:

\begin{definition}
    The \textit{separated Teichm\"{u}ller space} of $N$ is
    $$\Teich_s(N) = H(\Teich(N)),$$
    the Hausdorffification of the Teichm\"{u}ller space $\Teich(N)$. Concretely, any two inseparable points in $\Teich(N)$ are identified in the separated quotient $\Teich_s(N)$.
\end{definition}

The local Torelli theorem tells us that the induced period mapping 
$$\mathcal{P}_s: \Teich_s(N) \to \mathcal{D}$$
will still be a local isomorphism. The separated Teichm\"{u}ller space will be the space on which our first global Torelli theorem will be stated and shown. It is also useful to define two other moduli spaces as they will naturally emerge in our analysis of Enriques manifolds. 

\begin{definition}
    The \textit{K\"{a}hler Teichm\"{u}ller space} of $N$ is
    $$\Teich_{\HK}(N) = \{\text{marked K\"{a}hler structures on $N$}\}/\Diff_0(N),$$
    i.e. the moduli space of marked smooth isotopy classes of complex structures $[Y, \tilde{\varphi} : \tilde{N} \to (X,\tilde{\omega})]$, where $\tilde{\omega}$ is the image of a chosen K\"{a}hler class $\omega$ on the complex Enriques manifold $Y$ under the transfer map.
\end{definition}

\begin{remark}
    The subscript ``$\HK$" is used to mirror the notation used by Looijenga when he defines a \textit{hyperk\"{a}hler Teichm\"{u}ller space} \cite[pg. 9-10]{looijenga2021teichmuller}. Enriques manifolds certainly do not admit hyperk\"{a}hler structures, as they do not support any nondegenerate holomorphic $2$-forms.
\end{remark}

By Yau's theorem \cite{yau1977calabi}, every K\"{a}hler class $\omega$ on a complex Enriques manifold $Y$ contain a uniquely compatible Ricci-flat representative. This leads us to define a moduli of marked metric structures on $N$:

\begin{definition}
The \textit{metric Teichm\"{u}ller space} of $N$ is
$$\Teich_{\RF}(N) = \{\text{marked unit volume Ricci-flat structures on $N$}\}/ \Diff_0(N),$$
i.e. the moduli space of marked smooth isotopy classes of unit volume Einstein metric structures on $N$.
\end{definition}

By basic covering space theory, the universal cover of a Ricci-flat Enriques manifold $(N,g)$ is the Ricci-flat manifold $(\tilde{N}, \tilde{g})$, which is unique up to isometry. Using the twistor construction, we can assign to $(\tilde{N}, \tilde{g})$ its $D$-invariant period data
$$P_g := P_{\tilde{g}} \in \Gr_3^+ (H_\RR)$$
which, by \cite[Lemma 2.3]{oguiso2011enriques}, is endowed with a faithful orthogonal representation of the fundamental group 
$$D \hookrightarrow \SO(P_g) = \SO(3).$$
\begin{definition}
The \textit{K\"{a}hler period manifold} of an Enriques manifold $Y$ is the space
$$\mathcal{D}_{\HK}(Y) := \{(\omega, P_g) \in  \PP(H^{1,1}(Y,\RR))  \times \Gr_3^+(H_\RR) : \tilde{\omega} \in \mathcal{D}(P_g)\}.$$
\end{definition}

All of these Teichm\"{u}ller spaces come equipped with natural $\MOD(N)$-actions defined by pullback. In complete analogy with the hyperk\"{a}hler setting, we can build natural period domains and projection maps associated to each of these Teichm\"{u}ller spaces of Enriques manifolds by leveraging the twistor construction on hyperk\"{a}hler manifolds. These morphisms will also respect the $\MOD(N)$-actions. Following \cite[Corollary 2.7]{looijenga2021teichmuller}, we have that Property (3) of \Cref{TheMainLemma} and descent along the projection maps yields
\begin{proposition}\label{TeichProjections}
    $\Teich_{\HK}(N)$ is a separated complex manifold. The natural projection $\Teich_{\HK}(N) \to \Teich(N)$ is open with fibers open convex subsets of a hyperbolic space. The projection $\Teich_{\HK}(N) \to \Teich_{\RF}(N)$ is a $2$-sheeted cover.
\end{proposition}
\begin{proof}
    It is clear from the twistor construction that both maps are open submersions. The fibers of the projection $\Teich_{\HK}(N) \to \Teich(N)$ are given by the rays spanned by K\"{a}hler classes in the positive cone of $H^{1,1}(N,\RR)$. Hyperbolicity follows from the Hodge index theorem for the Beauville--Bogomolov--Fujiki form, and can even be seen by the computations done on the known examples in \Cref{EnriquesReview}.

    Given a Ricci-flat metric $g$ on $N$, its lift to the universal cover determines a positive $3$-plane $P_g$ equipped with a $q$-orthogonal $D$-action. Since a faithful $D$-representation into $\SO(P_g)$ admits a unique invariant line $\ell$ in $P_g$, the two points given by $\mathcal{D}(P_g) \cap \ell$ determine two antipodal $D$-invariant K\"{a}hler-Einstein structures on $\tilde{N}$ from which $g$ arose. 
\end{proof}
For each pair of Teichm\"{u}ller spaces, there is a naturally induced lifting map coming from the smooth universal cover $p : \tilde{N} \to N$. Concretely, we use the fact that $p$ is a local diffeomorphism to pull back any marked metric/complex/K\"{a}hler-Einstein structure on $N$ to the universal cover $\tilde{N}$. This produces the following commutative diagram
\begin{center}
\begin{tikzcd}
\Teich(N) \arrow[d, "p^*"']                       & \Teich_{\HK}(N) \arrow[d, "p^*"'] \arrow[r] \arrow[l]                       & \Teich_{\RF}(N) \arrow[d, "p^*"]                       \\
\Teich(\tilde{N}) \arrow[d, "\widetilde{\mathcal{P}}"'] & \Teich_{\HK}(\tilde{N}) \arrow[d, "\widetilde{\mathcal{P}}_{\HK}"'] \arrow[l] \arrow[r] & \Teich_{\RF}(\tilde{N}) \arrow[d, "\widetilde{\mathcal{P}}_{\RF}"] \\
\mathcal{D}(X)                              & \mathcal{D}_{\HK}(X) \arrow[l] \arrow[r]                                       & \Gr_{3}^{+}(H_{\RR})         
\end{tikzcd}
\end{center}

\begin{proposition}\label{invariantPeriodImage}
    The image of $\widetilde{\mathcal{P}}\circ p^*$ lies in $\mathcal{D}(Y)$, the image of $\widetilde{\mathcal{P}}_{\HK} \circ p^*$ lies in $\mathcal{D}_{\HK}(Y)$, and the finally the image of $\widetilde{\mathcal{P}}_{\RF}\circ p^*$ lies in $\Gr_3^+(H_\RR)^D$.
\end{proposition}
\begin{proof}
For a geometric structure on $\tilde{N}$ to have come from lifting, i.e. to lie in the image of $p^*$, it must be $D$-invariant in the correct sense, which we now describe. For a given $D$-invariant complex structure on $X$, we see that the nonvanishing holomorphic $2$-form on $X$ must lie in $H_{\CC,1}$. Given a $D$-invariant Ricci-flat metric on $\tilde{N}$, the period data on $(\tilde{N},\tilde{g})$ defined by the positive $3$-plane $P_g$ is fixed set-wise under the faithful $D$-action. These descriptions of $D$-invariance must have descended from the hyperk\"{a}hler period mapping.
\end{proof}

A key point in our proof of the smooth Birman--Hilden theorem for hyperk\"{a}hler manifolds, just as with the MacLachlan--Harvey proof \cite{maclachlan1975mapping}, is that the map of Teichm\"{u}ller spaces $p^*$ should be injective. In the classical case of Riemann surfaces, we see this by showing that Teichm\"{u}ller geodesics in $\Teich(Y)$ map to Teichm\"{u}ller geodesics in $\Teich(X)$ of the same length. In our original case of interest, e.g. when $X$ is a K3 surface and $Y$ is an Enriques surface, different technology is required to verify injectivity of the lifting map. We have laid the groundwork to show this in the more general setting of hyperk\"{a}hler manifolds. The following is the precise statement of \Cref{IntroThmB}:

\begin{theorem}[Global Torelli Theorem for Enriques manifolds]\label{Etorelli}
    The lifting map of Teichm\"{u}ller spaces $p^* : \Teich(N) \to \Teich(\tilde{N})$ is injective. As a consequence, the period mapping
    $$\mathcal{P}_{s}: \Teich_s(N) \to \mathcal{D}(Y)$$
    is injective on connected components of $\Teich_s(N)$. In
particular, the $\MOD(N)$-stabilizer of a component of $\Teich(N)$ acts with finite kernel on $H^2(N,\ZZ)$.
\end{theorem}
\begin{proof}
    Let us first show injectivity of the lifting map. If $p^*(x) = p^*(y)$, that means the $D$-invariant marked complex structures on $\tilde{N}$ are isotopic. By the proof of Property (1) in \Cref{TheMainLemma}, such $D$-invariant marked hyperk\"{a}hlerian manifolds lie in a $D$-invariant family over a connected base $B$, and so they can be made symmetrically isotopic. Thus the isotopy class of marked complex structure descend to $\Teich(N)$, showing us that $x= y \in \Teich(N)$. Thus $p^*$ is injective.
    
    For the second part, \Cref{invariantPeriodImage} tells us that period map $\mathcal{P}_{s}$ can be factored as
    $$\widetilde{\mathcal{P}}_s \circ p^* : \Teich(N) \to \Teich(\tilde{N}) \to \mathcal{D}(Y) \subset \mathcal{D}(X).$$
    The global Torelli theorem for hyperk\"{a}hlerian manifolds (\cite[Theorem 3.1]{looijenga2021teichmuller}) tells us that $\widetilde{\mathcal{P}}_s$ is injective on connected components. Since $\mathcal{P}_{s}$ is the composition of injective functions on connected components, we have proven injectivity of $\mathcal{P}_s$ on connected components. 
    
    Finally, given the stabilizer $\MOD(N)_\mathcal{C}$ of a connected component $\mathcal{C}\subset \Teich(N)$, we would like to show that it acts with finite kernel on $H^2(N,\ZZ)$. We follow the argument given in the proof of \cite[Claim 2.1]{verbitsky2019mapping}. Let $G < \MOD(N)_\mathcal{C}$ be the subgroup of mapping classes which act trivially on $H^2(N,\ZZ)$. Since $G$ commutes with the period mapping $\mathcal{C} \to \mathcal{D}(Y)$, and $G$ acts trivially on $\mathcal{D}(Y)$, it must act trivially on $\mathcal{C}$ itself. Thus we may pick a complex structure $J \in \mathcal{C}$ preserved by $G$, i.e. $G$ acts by holomorphic automorphisms on $(N,J)$, and moreover acts trivially on $H^{1,1}(N,J)$ by assumption. By the Calabi--Yau theorem, $G$ fixes a K\"{a}hler-Einstein metric, and since the group of isometries of a closed Riemannian manifold is compact, we can conclude that the group $G$ is finite.
\end{proof}

To the best of my knowledge, this is the first known instance of a global Torelli theorem for generalized Enriques manifolds to appear in the literature. Given the proof of Property (1) of \Cref{TheMainLemma}, and that we have descent along projection maps between Teichm\"{u}ller spaces, the lifting maps associated to the other Teichm\"{u}ller spaces are injective as well. More concretely, the choice of a compatible K\"{a}hler class on $Y$ lifts to a unique $D$-invariant K\"{a}hler class on $X$ by transfer, so the map of hyperk\"{a}hler Teichm\"{u}ller spaces must be injective. The twistor construction tells us that the Ricci-flat metric associated to any hyperk\"{a}hler structure must also then uniquely lift throught the map of Ricci-flat Teichm\"{u}ller spaces. To summarize:

\begin{corollary}
    The lifting maps of metric and k\"{a}hler Teichm\"{u}ller spaces $\Teich_{\RF}(N) \to \Teich_{\RF}(\tilde{N})$ and $\Teich_{\HK}(N) \to \Teich_{\HK}(\tilde{N})$ are injective.
\end{corollary}

The following result characterizes the geometry of the image of $\mathcal{P}_{\RF}$:

\begin{corollary}\label{geodesic}
The image of the lifting map $p^*(\Teich_{\RF}(Y)) \subset \Teich_{\RF}(X)$ is totally geodesic inside $\Teich_{\RF}(X)$ with respect to the symmetric metric $g_{loc}$ on each connected component.
\end{corollary}
\begin{proof}
Recall that the $\MOD(X)$-action on $\Teich_{\RF}(X)$ leaves the Beauville--Bogomolov--Fujiki form invariant \cite[Theorem 3.5]{verbitsky2013mapping}, which defines a symmetric metric $g_{loc}$ on each connected component of $\Teich_{\RF}(X)$ through the period mapping $\widetilde{\mathcal{P}}_{\RF}$. Each connected component is isometric to one another, and looks like their image in $\Gr_3^+(H_\RR)$. The restriction of the $\MOD(X)$-action to the $D$-action on $\Teich_{\RF}(X)$ is thus isometric with respect to the symmetric metric $g_{loc}$ on $\Teich_{\RF}(X)$. Recall that the connected components of the fixed set of an isometric action on a Riemannian manifold form a totally geodesic submanifold \cite[Chapter II, Theorem 5.1]{kobayashi1972transformation}. The image of $p^*$ agrees with the fixed set of the isometric $D$-action on $\Teich_{\RF}(X)$, the $D$-invariant metrics on $X$, and thus is totally geodesic. 
\end{proof}

Injectivity of the lifting map of Teichm\"{u}ller spaces also implies the existence of almost-universal families of Enriques manifolds. The following results for Enriques manifolds are analogous to the results of Markman for hyperk\"{a}hler manifolds \cite{markman2021existence}. We model our statements off of \cite[Corollary 4.5, 4.6]{looijenga2021teichmuller}:

\begin{theorem}
    The Teichm\"{u}ller space of Einstein metrics on $N$, $\Teich_{\RF}(N)$, carries a family of Einstein manifolds $\mathcal{U}_{\RF}(N) /\Teich_{\RF}(N)$ which is endowed with a faithful action of $\MOD(N)$. It is almost-universal in the sense that every family of Einstein metrics on $N$ is a pull-back of this one, but can be so in more ways than one, with the ambiguity residing in a finite group which is constant on every connected component of $\Teich_{\RF}(N)$.
\end{theorem}
\begin{proof}
If we pull back the universal family of Einstein manifolds diffeomorphic to $\tilde{N}$ along the injective map $p^*$: 
    \begin{center}
        \begin{tikzcd}
p^* \mathcal{U}_{\RF}(\tilde{N}) \arrow[d] \arrow[r] & \mathcal{U}_{\RF}(\tilde{N}) \arrow[d] \\
\Teich_{\RF}(N) \arrow[r, "p^*"']                    & \Teich_{\RF}(\tilde{N})               
\end{tikzcd}
    \end{center}
we obtain an almost-universal family of $D$-invariant Einstein manifolds over $\Teich_{\RF}(N)$. By construction, their $D$-quotients are Ricci-flat Enriques manifolds diffeomorphic to $N$. Taking the $D$-quotient of the family $p^* \mathcal{U}_{\RF}(\tilde{N})$ fiberwise produces our desired almost-universal family $\mathcal{U}_{\RF}(N) /\Teich_{\RF}(N)$.
\end{proof}

To quote Looijenga \cite[pg. 15]{looijenga2021teichmuller}, ``in somewhat fancier language: $\Teich_{\RF}$ underlies a (Deligne-Mumford) stack and this
stack is a constant gerbe on every connected component." By \Cref{TeichProjections}, fibers of the natural projection maps out of $\Teich_{\HK}$ come with tautological families, which implies:

\begin{corollary}
    The Teichm\"{u}ller spaces $\Teich_{\HK}(N)$ and $\Teich(N)$ support families above them that are almost-universal in the sense described above.
\end{corollary}

\begin{question}
    What are the images of these period mappings? We refer to \cite{namikawa1985periods} for the case of Enriques surfaces; it is not entirely obvious, even given what was known for K3 surfaces. 
    
    The images of period maps for hyperk\"{a}hler manifolds are described in \cite[Theorem 4.9]{amerik2015teichmuller} and \cite[Section 4]{looijenga2021teichmuller}. This tells us, for example, that the image of $\widetilde{\mathcal{P}}_{\RF}$ lies in the complements of divisors determined by MBM classes. What else is needed to characterize the image for Enriques manifolds?
\end{question}
\section{The Birman--Hilden Type Isomorphism}\label{BHIsomorphism}

We shall now record some key properties of the $\MOD(M)$-actions on Teichm\"{u}ller spaces. Understanding these actions will naturally lead to our proof of the Birman--Hilden theorem.
\begin{proposition}\label{Faithful}
Let $M$ be the smooth manifold underlying either a hyperk\"{a}hler or Enriques manifold. The mapping class group $\MOD(M)$ acts faithfully on the Teichm\"{u}ller space $\Teich_{\RF}(M)$.
\end{proposition}
\begin{proof}
Assume for the sake of nontriviality that $f \in \MOD(M)$ preserves a connected component $\mathcal{C} \subset \Teich_{\RF}(M)$. Property (4) of \Cref{TheMainLemma} for Enriques manifolds, and Property (4) of \cite[Proposition 2.3]{looijenga2021teichmuller} for hyperk\"{a}hler manifolds, imply that 
$$G := \ker(\MOD(M)_\mathcal{C}\to  \Aut(H^2(M,\ZZ)))$$
can be identified with the finite group of isotopically trivial isometries $\Isom_0(M,h)$ for any $h \in \mathcal{C}$. We know that $G$ must lie in the Torelli group, but if $G$ is given by a group of isotopically trivial isometries, its image in the mapping class group must be trivial. Thus $\MOD(M)$ acts faithfully on $\Teich_{\RF}(M)$
\end{proof}

\begin{lemma}\label{SymmetricAction}
Let $Y$ be a smooth Enriques manifold with universal cover $X$. The symmetric mapping class group $\SMod(X)$ acts on the space of $D$-invariant Ricci-flat hyperk\"{a}hler manifolds $p^*(\Teich_{\RF}(Y))$ with kernel given by the isotopy group of deck group elements.
\end{lemma}
\begin{proof}
Since every element of $\SMod(X)$ is fiber-preserving, the action of $\SMod(X)$ on $p^*(\Teich_{\RF}(Y))$ is well-defined because the pullback of a $D$-invariant metric will remain $D$-invariant. The deck group $D$ is equal to the kernel for the following reason. Suppose $\alpha \in \SMod(X)$ acts by the identity on $p^*(\Teich_{\RF}(Y))$. Pick a representative symmetric diffeomorphism $\tilde{f} \in \alpha$ and consider its pushforward $f \in \Diff(Y)$. \Cref{Etorelli} allows us to identify $\Teich_{\RF}(Y)$ with $p^*(\Teich_{\RF}(Y))$, telling us that $f$ acts by the identity on $\Teich_{\RF}(Y)$. It follows that $f$ is isotopic in $\Diff(Y)$ to the identity diffeomorphism, by faithfulness of the $\MOD(Y)$-action shown in \Cref{Faithful}. By \Cref{symlift}, an isotopy between $f$ and $\id$ lifts, once we have selected basepoints, to a symmetric isotopy between $\tilde{f}$ and $\tilde{\id}$. We can pull back this isotopy by any one of the elements of the deck group $\eta \in D$ to get an isotopy from a lift of $f$ to $\eta$. This proves the claim.
\end{proof}

We can now state and prove the main result, which is \Cref{IntroThmA} in the introduction.

\begin{theorem}
Suppose we have a smooth regular cover $p: X \to Y$ with $X$ a hyperk\"{a}hlerian manifold, generated by a deck group $D$. Then the cover $p$ has the Birman--Hilden property. That is, isotopic symmetric diffeomorphisms on $M$ can be made symmetrically isotopic with respect to the $D$-action.
\end{theorem}
\begin{proof}
By \Cref{SymmetricAction}, the symmetric mapping class group $\SMod(X)$ acts on $p^*(\Teich_{\RF}(Y))$ by restricting the $\MOD(X)$-action on $\Teich_{\RF}(X)$, with kernel equal to $D$. By \Cref{Faithful}, the mapping class group $\MOD(X)$ acts faithfully on $\Teich_{\RF}(X)$. Moreover, injectivity of the lifting map $p^*$ (\Cref{Etorelli}) tells us that $\MOD(Y)$ also acts faithfully on $p^*(\Teich_{\RF}(Y))$, and the action must factor through the $\SMod(X)$-action. Since $\MOD(Y) \to \SMod(X)/D$ is injective, it follows from \Cref{BHequivalence} that $\SMod(X)/D \cong \MOD(Y)$, so $p:X \to Y$ has the Birman--Hilden property, thereby proving the claim.
\end{proof}

\begin{question}[Relative Mapping Class Groups]
    The unbranched version of the Birman-Hilden theorem in the case of Riemann surfaces is quite simple; dealing with branch loci makes the problem more complicated. Recall that every K3 surface $M$ is diffeomorphic to a double cover of $\CC\PP^2$ branched over a smooth sextic curve $C \subset \CC\PP^2$. Is there a branched version of the Birman--Hilden theorem in this setting, something like an isomorphism $\SMod(M) \cong \MOD(\CC\PP^2,C)$? As far as I am aware, the techniques used in this paper do not immediately apply, as there isn't an obvious geometry one can leverage on a ``Teichm\"{u}ller space of sextic curves in the plane". See \cite[pg. 25]{hain2023mapping} for related problems.
\end{question}

\section{Nielsen Realization on Enriques Surfaces}\label{NielsenSection}

We now will apply the hyperk\"{a}hler Birman--Hilden theorem to Nielsen realization problems for Enriques surfaces. Fix a smooth manifold $M$. The Nielsen realization problems ask the following questions:

\begin{question}[The Metric Nielsen Realization Problem]
Let $G \subset \MOD(M)$ be a subgroup and let $\mathcal{M}$ be a distinguished class of metrics on $M$. Does there exist a lift of $G \subset \Diff(M)$ and a metric $g \in \mathcal{M}$ on $M$ such that $G \subset \Isom(M,g)$?
\end{question}

\begin{question}[The Complex Nielsen Realization Problem]
Let $G \subset \MOD(M)$ be a subgroup. Does there exist a lift of $G \subset \Diff(M)$ and a complex structure $J$ on $M$ such that $G \subset \Aut(M,J)$?
\end{question}

To begin answering these questions, we recount some specific details regarding Enriques surfaces and their covering K3 surfaces below. A good reference is \cite{huybrechts2016lectures}. Let $M$ be the smooth $4$-manifold which underlies a K3 surface. By Poincar\'{e} duality and naturality of the cup product, $M$ has a canonical orientation $[M] \in H_4(M,\ZZ)$ which is preserved by any self-diffeomorphism of $M$. Let us equip $H = H^2(M,\ZZ)$ with the quadratic form
$$q(\alpha) := \frac{1}{2} \int_{[M]} \alpha \cup \alpha.$$
The quadratic form $q$ then defines an even unimodular symmetric bilinear form of signature $(3,19)$, and in fact agrees with the Beauville--Bogomolov--Fujiki form in this low-dimensional setting. This turns $H$ into the unique such lattice of signature $(3,19)$, the \textit{K3 lattice}:
$$(H,q) \cong U\oplus U\oplus U\oplus E_8(-1) \oplus E_8(-1).$$
This isomorphism of lattices is not canonical. Let $p : M\to N$ be a universal covering of the Enriques manifold $N$ by the K3 manifold $M$ generated by the involutive deck group $D = \langle \iota \rangle$. There exists some lattice basis on $H^2(M,\ZZ)$ for which the induced action of $\iota$ on $H$ is given by
\begin{align*}
    \iota^* :U\oplus U\oplus U\oplus E_8(-1) \oplus E_8(-1) &\longrightarrow U\oplus U\oplus U\oplus E_8(-1) \oplus E_8(-1)\\
    (x_1,x_2,x_3,x_4,x_5) &\longmapsto(-x_1,x_3,x_2,x_5,x_4).
\end{align*}
The transfer homomorphism
$$p^* : H^2(N,\ZZ) \to  H^2(M,\ZZ)$$
is injective on the free part of cohomology $H^2(X,\ZZ)_{\mathrm{free}} := H^2(X,\ZZ) / \mathrm{torsion}$ with image isomorphic to the invariant sublattice
$$H^{\langle\iota\rangle} \cong U(2) \oplus E_8(-2).$$
Thus $H^2(M,\ZZ)_{\mathrm{free}}$ is isomorphic to $U \oplus E_8(-1)$, which we call this the \textit{Enriques lattice}. It is the unique such even unimodular lattice of signature $(1,9)$. 

Let us now recall some results regarding Nielsen realization for K3 surfaces $M$ due to Farb--Looijenga \cite{farb2021nielsen}. Let $G \subset \OO(3,19)^+(\ZZ)$ be a finite subgroup. Since $\OO(3,19)^+(\RR)$ has the nonpositively curved Grassmannian $\Gr(3,19)(\RR)$ as its associated symmetric space, the natural $G$-action on $\Gr(3,19)(\RR)$ has a fixed point, i.e. $G$ leaves invariant a positive $3$-plane $P$ in $H^2(M,\RR)$. It also preserves the orientation on $P$, thus defining a representation $\rho : G \to \SO(P)$. Such representations can be cyclic, dihedral, tetrahedral, octahedral, or icosahedral. The last three cases are irreducible. From this, we construct a representation-theoretic invariant of $G$ on $H^2(M,\RR)$:
\begin{definition}
Let $I_G \subset H^2(M,\RR)$ be the direct sum of all irreducible $\RR[G]$-submodules of the type that appear in $P$. This subspace $I_G$ is canonically associated to $G$. The \textit{realization invariant $L_G$} is defined to be
$$L_G := I_G^\perp \cap H^2(M,\ZZ).$$
It is clearly negative definite and $G$-invariant.
\end{definition}

The metric/complex realization result in \cite{farb2021nielsen} is as follows:
\begin{theorem}[K3 Metric/Complex Nielsen Realization]
Let $G \subset \OO(3,19)^+(\ZZ)$ be a finite group and $M$ the K3 manifold. The following holds true:
\begin{enumerate}
    \item $G$ lifts to the group of isometries of a Ricci-flat metric on $M$ if and only if $L_G$ contains no $(-2)$-vectors.
    \item $G$ lifts to the group of automorphisms of some complex structure on $M$ endowed with a Ricci-flat K\"{a}hler metric (making it a K\"{a}hler-Einstein K3 surface) if and only if $L_G$ contains no $(-2)$-vectors and in addition the trivial representation appears in $L_G^\perp$. In this case the complex structure can be chosen so that M is projective and $G$ acts by algebraic automorphisms.
\end{enumerate}
If instead $G$ is given rather as a finite subgroup of $\MOD(M)$, then a necessary and sufficient condition for it to lift as in Case 1 (resp. Case 2) is that $G$ satisfies the homological conditions stated in that case and also preserves a connected component of the Teichm\"{u}ller space
$\Teich_{\RF}(M)$.
\end{theorem}

This theorem, paired with the hyperk\"{a}hler Birman--Hilden theorem, allows us to pin down sufficient and necessary conditions for finite order elements of the mapping class group to be realizable by metric/complex automorphisms on the Enriques surface $N$, as given in \Cref{IntroThmC}. Thus metric/complex realizability in $\MOD(N)$ can be determined by linear algebra:

\begin{theorem}[Enriques Metric/Complex Cyclic Nielsen Realization]\label{EnriquesNielsen}
Let $N$ be a smooth Enriques surface with covering K3 surface $M$ and associated deck group $D$. Let $G \subset \MOD(N)$ be a finite subgroup of the mapping class group. Then $G$ lifts to a group of Ricci-flat isometries (respectively K\"{a}hler-Einstein automorphisms) if and only if its group of lifts $\tilde{G} \subset \SMod(M)$ is Ricci-flat (respectively K\"{a}hler-Einstein) realizable on $M$.
\end{theorem}
\begin{proof}
$(\Longrightarrow)$ This is clearly true, because any lift $\tilde{f} \in \tilde{G}$ of $f \in G$ will preserve the lifted metric/complex structure.\\
$(\Longleftarrow)$ The hyperk\"{a}hler Birman--Hilden theorem tells us that $\SMod(M)$ is a central extension of $\MOD(N)$ by $D$, i.e there is a short exact sequence 
$$1 \to D \to \SMod(M) \xrightarrow[]{\Phi} \MOD(N) \to 1.$$
This sequence restricts along the group $G$ to the group of lifts $\tilde{G}$ as
$$1 \to D \to \tilde{G} \to G \to 1.$$
The group $\tilde{G}$ consists of symmetric mapping classes, so any Ricci-flat (respectively K\"{a}hler-Einstein) structure that is preserved by $\tilde{G}$ is $D$-invariant, and thus will descend to a Ricci-flat (respectively K\"{a}hler-Einstein) structure on the Enriques surface $N$. Pushing $\tilde{G}$ forward through $\Phi$ gives us that $G$ preserves the Ricci-flat (respectively K\"{a}hler-Einstein) structure on $N$. To see this, note that if there was an element $f \in G$ that did not act by automorphisms, it would have a lift $\tilde{f}$ that would fail to preserve the lifted $D$-invariant geometric structure we started with, a contradiction. This proves the claim.
\end{proof}

We can now show that Einstein metric realizability implies complex realizability in the setting of Enriques surfaces, as mentioned in \Cref{IntroThmE}. This is complementary to the situation for K3 surfaces \cite[Theorem 1.4]{farb2021nielsen}:

\begin{theorem}[Ricci-flat implies complex]
Suppose a finite subgroup $G \subset \MOD(N)$ can be realized by isometries of some
Ricci-flat metric on the Enriques manifold $N$. Then it must also preserve a compatible K\"{a}hler-Einstein structure on $N$.
\end{theorem}
\begin{proof}
Given that $G$ is realizable within the group of isometries of some Ricci-flat metric $h$ on $N$, we will identify $G$ with its image in the isometry group $G \subset \Isom(N,h)$. \Cref{EnriquesNielsen} implies that $G$ lifts to its group of symmetric isometries $\tilde{G}$ of the Ricci-flat K3 surface $(M,\tilde{h})$, so $\tilde{G} \subset \Isom(M,\tilde{h})$. As such, $\tilde{G}$ must preserve the period data $P_{\tilde{h}} \in \Gr_3^+(H_\RR)$. This gives us a representation
$$\tilde{G} \to \SO(P_{\Tilde{h}}) = \SO(3).$$
Since $\tilde{G}$ consists of symmetric isometries, any $\tilde{f} \in \tilde{G}$ and the Enriques involution $\iota$ commute. It follows that their images in $\SO(3)$ must commute. Since $\rank_{\RR}(\SO(3)) = 1$, the two elements $\Tilde{f}^*$ and $\iota^*$ must share a common $1$-eigenspace. There is a unique such line $\ell\subset P_{\tilde{h}}$ determined by $\iota$, which by the twistor construction yields a compatible $D$-invariant K\"{a}hler-Einstein structure on $M$ preserved by $\tilde{G}$. By $D$-invariance, this then descends to a K\"{a}hler-Einstein structure on $N$ that $G$ preserves.
\end{proof}

We shall now study the problem of realization for order 2 subgroups of $\MOD(X)$ by smooth involutions. As with the previous work, we can leverage the results found in \cite{farb2021nielsen} when restricting our attention to the symmetric mapping classes on the K3 surface $M$. 

\begin{theorem}[Smooth involutive realization for Enriques surfaces]
Let $f \in \MOD(N)$
be an order 2 element of the Enriques mapping class group. Then $f$ is realized by a smooth involution of $N$ if and only if a lift into $\SMod(M)$ is smoothly realizable.
\end{theorem}
\begin{proof}
$(\Longrightarrow)$ This is clearly true, because we can smoothly lift $f$ through the universal covering map.\\
$(\Longleftarrow)$ By the hyperk\"{a}hler Birman--Hilden theorem, $f \in \MOD(N) \cong \SMod(M)/D$, so we can lift $G = \langle f\rangle $ into its group of lifts $\tilde{G} \subset \SMod(M)$. The group $\tilde{G}$ consists of symmetric mapping classes upstairs, so a smooth realization of $\tilde{G}$ will be $D$-invariant, and hence descend to a smooth realization on the Enriques surface $N$. This proves the claim.
\end{proof}
Another fundamental example of $4$-manifold diffeomorphisms we are interested in studying are Dehn twists about $(-2)$-spheres. That is, twisting via the geodesic flow around an embedded sphere $i : S^2\hookrightarrow M$ with $v=i_\ast[S^2]$ such that $v^2=-2$.

It is well known that $T^2_{S^2}$ is smoothly isotopic to the identity, which one can see by thinking about the nontrivial element of $\pi_1(\mathrm{SO}(3))$. Furthermore, $T_{S^2}$ induces a reflection across a hyperplane in cohomology. To elaborate, let $v=[S^2]\in H_2(M;\ZZ)$. Within $H^2(M;\RR)$, $(T_{S^2})_\ast$ preserves $v^\perp$ and flips $v\mapsto -v$. From our smooth Birman--Hilden theorem, we are able to conclude \Cref{IntroThmD}:

\begin{corollary}
Let $N$ be an Enriques surface with a universal K3 covering $p : M \to N$. Let $T_S \in \Diff(N)$ be the Dehn twist about an embedded $2$-sphere $S \subset N$ with self intersection $S \cdot S = -2$. Although $T_S^2$ is smoothly isotopic to the identity, $T_S$ is not topologically isotopic to any finite order diffeomorphism of $N$.
\end{corollary}
\begin{proof}
Suppose that $T_S$ sits in the same isotopy class of some finite order diffeomorphism $f : N \to N$ with $f^* = (T_S)^*$. Then $f$ acts by reflection about the $(-2)$-vector corresponding to $S$ on cohomology. 

An isotopy between $f$ and a representative twist $T_S$ will then lift to a symmetric isotopy between $\tilde{f}$ and the diffeomorphism $\tilde{T}_{S_1 + S_2}$, where the sum of the embedded spheres $S_1$ and $S_2$ forms the image of $S$ under $p^* :H^2(N,\ZZ) \to H^2(M,\ZZ)$. The action on cohomology is given by two reflections across distinct $(-2)$-vectors by transfer. Thus the finite even order diffeomorphism $\tilde{f}$ is a realization of the finite order 2 mapping class $\tilde{T}_{S_1 + S_2}$ upstairs, contradicting \cite[Proposition 3.7]{farb2021nielsen}: realizable involutive diffeomorphisms of K3 surfaces cannot have cyclotomic summands in cohomology. Promoting this claim to exclude a topological isotopy of $\tilde{T}_{S_1+S_2}$ to a finite order diffeomorphism follows from \cite[Theorem 1.12]{farb2021nielsen}
\end{proof}

\begin{question}[Lifting Isotopy Group Actions]
    Our results on Nielsen realization for Enriques surfaces $Y$ use a lifting argument coming from the hyperk\"{a}hler Birman--Hilden theorem, which tells us that $\SMod(X)$ is a central extension of $\MOD(Y)$ by $D = \langle \iota \rangle$. One can readily lift finite order mapping classes to the universal cover $X$ and study Nielsen realization element-wise, but it is a much more subtle problem to lift a whole $G$-action on $Y$ to an $\langle \iota \rangle$-equivariant $G$-action on $X$. Can one prove (or more likely disprove) such lifts of group actions exist? The best we could hope for (which perhaps is not true) is to show the group cohomology $H^2(\MOD(Y),\ZZ/2\ZZ)$ is trivial, i.e. $\SMod(X)$ is a direct product and the full mapping class group $\MOD(Y)$ trivially lifts.
\end{question}

\section{A Conjecture on the Enriques Surface's Torelli Group}\label{TorelliConjecture}

The following conjecture was suggested to me by Looijenga \cite{looijenga2023talking}:

\begin{conjecture}[Looijenga]
Let $\mathcal{I}_Y := \ker(\MOD(Y) \to \Aut(H^*(Y,\ZZ),Q_Y))$ be the smooth Torelli group of an Enriques surface $Y$. Then $\mathcal{I}_Y = \{1\}$, that is, the smooth Torelli group is trivial.
\end{conjecture}

Let us frame this conjecture within the context of exotic structures on $4$-manifolds. It is now well known that there exist infinitely many distinct homotopy K3 surfaces and homotopy Enriques surfaces, the latter by work of Okonek \cite{okonek1988fake}. Freedman's topological classification of $4$-manifolds \cite{freedman1982topology}, and its continuation carried out by Hambleton-Kreck \cite{hambleton1988classification}, implies that such manifolds are homeomorphic but not diffeomorphic to the smooth manifolds underlying K3 surfaces and Enriques surfaces, respectively. On the other hand, the Torelli theorem for K3/Enriques surfaces (\cite{looijenga1980torelli}, \cite{pyatetskii1971torelli}, \cite{burns1975torelli}, and \cite{namikawa1985periods}) imply that any two smooth manifolds underlying a complex K3/Enriques surface must be diffeomorphic. Such results indicate that it is quite subtle and not obvious how to determine whether cohomologically trivial mapping classes are smoothly isotopic to the identity for K3/Enriques surfaces.

Cohomologically trivial complex automorphisms of Enriques surfaces were classified by Mukai--Namikawa \cite{mukai1984automorphisms}. One can check directly that they are smoothly isotopic to the identity. If one can show that all smooth Torelli elements can be deformed to one of the only two holomorphic examples, the conjecture follows. 

The Torelli theorem for K3 surfaces implies that no complex structure can be preserved by any nontrivial Torelli mapping class, since the smooth Torelli group acts simply-transitively on $\pi_0(\Teich(X))$ \cite[Remark 3.5]{looijenga2021teichmuller}. Simple-transitivity of the $\mathcal{I}_Y$-action on $\pi_0(\Teich(Y))$ also must hold, since the moduli space of Enriques surfaces is connected \cite[Corollary 1.16]{namikawa1985periods}. Therefore we have the following:
\begin{theorem}\label{iffTorelli}
The smooth Torelli group of the Enriques manifold is trivial if and only if the smooth Torelli group of the K3 manifold is trivial.
\end{theorem}
\begin{proof}
As Looijenga remarks in \cite[pg. 12]{looijenga2021teichmuller}, establishing connectivity of the Teichm\"{u}ller space is equivalent to proving the triviality of the Torelli group $\mathcal{I}_X$ for K3 surfaces $X$. By \Cref{geodesic} and simple-transitivity of the Torelli action on $\pi_0(\Teich_{\RF})$, every connected component of $\Teich_{\RF}(X)$ contains a unique connected component of $\Teich_{\RF}(Y)$, as it is a totally geodesic subspace determined by the fixed set of the $D$-action on $\Teich_{\RF}(X)$. Thus one marked moduli space is connected if and only if the other is.
\end{proof}

The tools needed to determine whether there are nontrivial elements in the smooth Torelli group do not seem to be available yet. In fact, it is not clear to me whether we know the answer in the topological category either. For the K3 manifold $X$, theorems of Freedman \cite{freedman1982topology}, Kreck \cite{kreck2006isotopy}, Perron \cite{perron1984pseudo}, and Quinn \cite{quinn1986isotopy} imply that the natural representation 
$$\pi_0(\Homeo^+(X)) \to \Aut(H^2(X,\ZZ),Q_X) \cong \OO(3,19)(\ZZ)$$
is an isomorphism, so the continuous K3 Torelli group is trivial. Hambleton--Kreck proved that the analogous representation for the continuous mapping class group of the Enriques manifold $Y = X / \langle \iota \rangle$ is surjective \cite{hambleton1988classification}. I cannot find any further work that determines what lies in the kernel of
$$\pi_0(\Homeo^+(Y)) \to \Aut(H^2(Y,\ZZ), Q_Y).$$
A Birman--Hilden type theorem should be true in for these topological $4$-manifolds, and a proof is probably amenable to topological surgery methods.

%
%

\bibliographystyle{alpha} 
\bibliography{HKBirmanHilden}

\Addresses

\end{document}